\documentclass[11pt]{article}

\usepackage{geometry}   
\usepackage{amsmath}
\usepackage{amssymb}
\usepackage{amsthm}

\usepackage{hyperref}
\hypersetup{colorlinks=true,linkcolor=blue,pdfstartview=FitH}

\usepackage{graphicx} 
\usepackage{tikz}
\usetikzlibrary{matrix,arrows,positioning}

\newcounter{remark}
\newcommand{\remark}{\addtocounter{remark}{1}
                       \par \quad {\bf \arabic{remark}}.\,
                      }

\newenvironment{commentaire}{\begin{quote}
                     \normalfont\footnotesize
                     \setlength{\parskip}{0.25\parskip}
                     \renewcommand{\item}{\remark}
                     {\bf }
                     \nobreak
                    }
                    {\end{quote}}

\newenvironment{rk}{\begin{quote}
                     \normalfont\footnotesize {{\bf Remark} --}
                    }{\end{quote}}

\newenvironment{preuve}{\medbreak \noindent {\bf Proof~---}}
                       {\hfill $\square$ \medbreak}

\newenvironment{psmallmatrix}{\left(\begin{smallmatrix}}{\end{smallmatrix}\right)}
\newenvironment{vsmallmatrix}{\left|\begin{smallmatrix}}{\end{smallmatrix}\right|}

\newcommand{\NN}{\mathbb N}
\newcommand{\ZZ}{\mathbb Z}

\newcommand{\QQ}{\mathbb Q}
\newcommand{\RR}{\mathbb R}

\newcommand{\FF}{\mathbb F}
\newcommand{\Cun}{X}
\newcommand{\Cdeux}{Y}

\newcommand{\Dcal}{\mathcal W}

\newcommand{\Hcal}{\mathcal H}
\newcommand{\Lcal}{\mathcal L}

\newcommand{\Ocal}{\mathcal O}

\newtheorem{theo}{Theorem}

\newtheorem{defi}[theo]{Definition}
\newtheorem{cor}[theo]{Corollary}
\newtheorem{lem}[theo]{Lemma}

\newtheorem{prop}[theo]{Proposition}

\newtheorem{criteria}[theo]{Criteria}

\newcommand{\Det}{\operatorname{Det}}

\newcommand{\Gram}{\operatorname{Gram}}

\newcommand{\Id}{\operatorname{Id}}

\newcommand{\NS}{\operatorname{NS}}

\newcommand{\Vect}{\operatorname{Vect}}

\newcommand{\HV}{\Vect(H, V)}
\newcommand{\HVperp}{{\mathcal E}}
\newcommand{\Frob}{{\mathcal F}}
\newcommand{\gammaa}{\gamma_{12}}

\newcommand{\scalaire}[2]{\left\langle #1, #2 \right\rangle}

\title{From Hodge Index Theorem to the number of points of curves over finite fields}
\author{
Emmanuel Hallouin \& Marc Perret\thanks{Institut de Math\'ematiques de Toulouse, UMR 5219, hallouin@univ-tlse2.fr, perret@univ-tlse2.fr}
}
\date{\today}

\begin{document}

\maketitle

\begin{abstract} 

We push further the classical proof of Weil upper bound for the number of rational points of an absolutely irreducible smooth projective curve $X$ over a finite field in term of euclidean relationships between the Neron Severi  classes in $X\times X$ of the graphs of iterations of
the Frobenius morphism. This allows us to recover Ihara's bound, which can be seen as a {\em second order} Weil upper bound, to establish a new {\em third order} Weil upper bound, and using {\tt magma} to produce numerical tables for {\em higher order} Weil upper bounds. We also give some interpretation for the defect of exact recursive towers, and give several new bounds for points of curves in relative situation $X \rightarrow Y$.
\end{abstract}




\section*{Introduction}

Let $\Cun$ be an absolutely irreducible smooth projective curve defined over the finite field~${\mathbb F}_q$ with $q$ elements. The classical proof of Weil  Theorem for the number of $\FF_q$-rational points $\sharp \Cun({\mathbb F}_q)$ rests upon Castelnuovo identity \cite{Weil}, a corollary of Hodge index Theorem for  the smooth algebraic surface $\Cun \times \Cun$. The intent of this article  is to push further this viewpoint by forgetting Castelnuovo Theorem. We come back to the consequence of Hodge index
Theorem that the intersection pairing on the Neron Severi
space~$\NS(X\times X)_{\mathbb R}$ is anti-euclidean on the orthogonal complement~$\HVperp_{\Cun}$ of the trivial plane generated by the horizontal and vertical classes.
Thus, the opposite~$\scalaire{\cdot}{\cdot}$ of the intersection
pairing endows~$\HVperp_{\Cun}$ with a structure of euclidean space.
 Section~\ref{s_Euclidean_space} is devoted to
 few useful scalar products computations.

\medskip

In section~\ref{s_Absolute_bounds},
we begin by giving a proof of Weil inequality which is, although equivalent in principle, different in presentation than the usual one using Castelnuovo Theorem given for instance in \cite[exercice~1.9,~1.10 p. 368]{Hartshorne}
or in \cite[exercises 8, 9, 10 p. 251]{Shafarevich}). We prove that Weil bound follows from Cauchy-Schwartz inequality for the orthogonal projections
onto~$\HVperp_{\Cun}$ of the diagonal class~$\Gamma^0$ and the class~$\Gamma^1$ of the graph of the Frobenius morphism on~$\Cun$.


The benefit of using Cauchy-Schwartz instead of Castelnuovo is the following. We do not know what can be a Castelnuovo identity for more than two Neron Severi classes, while we do know what is Cauchy-Schwartz for any number of vectors. It is well-known that Cauchy-Schwartz between two vectors is the non-negativity of their $2\times 2$ Gram determinant.
 Hence, we are cheered on investigating the consequences of the non-negativity of larger Gram determinants involving the Neron Severi classes~$\Gamma^{k}$ of the graphs of the $k$-th iterations  of the Frobenius morphism. 


For the family~$\Gamma^0,\Gamma^1,\Gamma^2$,
we recover the well known Ihara bound~\cite{Ihara}
which improves Weil  bound for curves of genus greater than~$g_2=\frac{\sqrt{q}(\sqrt{q} - 1)}{2}$,
a constant appearing very
naturally with this viewpoint in section~\ref{s_Ihara_bound}, especialy looking at figure~\ref{figure_Ihara_domain} in section~\ref{s_Ihara_bound}. It follows that
the classical Weil  bound can be seen as a {\em first order Weil bound}, in that it comes from the
euclidean constraints between $\Gamma^0$ and $\Gamma^1$, while
the Ihara bound can be seen as a {\em second order Weil bound}, in that it comes from
euclidean constraints between $\Gamma^0$, $\Gamma^1$ and $\Gamma^2$.
Moreover this process can be pushed further: by considering the
family~$\Gamma^0,\Gamma^{1}, \Gamma^{2}$ and~$\Gamma^{3}$, we obtain a new
{\em third order Weil bound} (Theorem~\ref{Ordre3}), which improves the Ihara
bound for curves of genus greater than another
constant~$g_3= \frac{\sqrt{q}(q-1)}{\sqrt{2}}$.

The more the genus increases, the more the Weil bound should be chosen of high order
in order to be optimal. Therefore it is useful to compute higher order
Weil bounds. Unfortunately, establishing them explicitely requires the resolution of high degree one variable polynomial equations over ${\mathbb R}$. For instance, the usual Weil bound requires the resolution of a degree one equation, Ihara and the third order Weil bounds require the resolution of second order equations, while fourth and fifth order Weil bounds require the resolution of third degree equations, and so on. Moreover, a glance at Ihara {\em second order Weil bound} and at our explicit {\em third order Weil bound} will convince the reader that they become more and more ugly as the order increase. 

Hence, we give up the hope to establish explicit formulae for the Weil bounds of order from 4. 
We then turn to a more algorithmic point of
view. We use an algorithm which, for a given genus~$g$ and a given field
size~$q$, returns the best upper order Weil bound for the number of
${\FF}_q$-rational points of a
genus~$g$ curve, together with the corresponding best order $n$.
The validity of this algorithm requires some results proved in the second
section. In the figure below, we represent the
successive Weil bounds (in logarithmic scales) of order from~$1$ to~$5$. Note that taking into account the logarithmic scale for the $y$-axis, higher order Weil bounds become significantly better than usual Weil's one.

\begin{center}   
\includegraphics[scale=0.75]{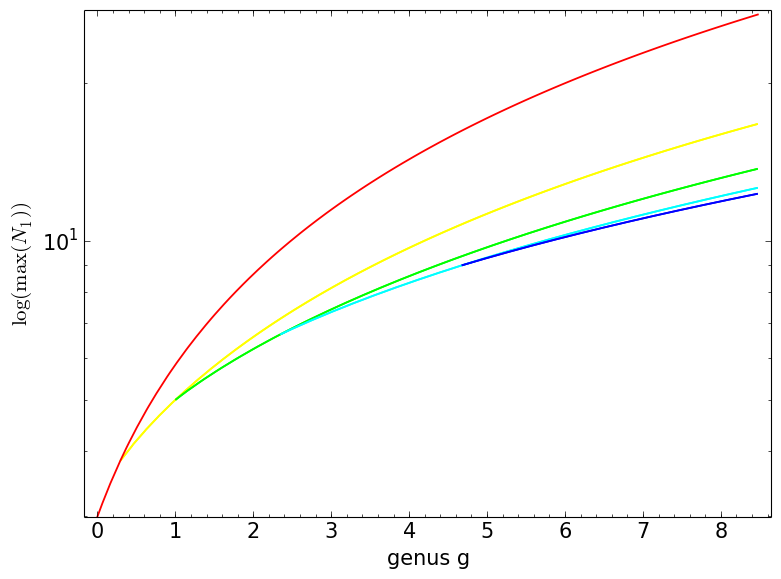}
\end{center}
\begin{commentaire} Weil bounds of order~$1$ to~$5$ for~$\sharp X(\FF_q)$ (here for $q=2$). Note that the $y$ axis is logarithmic. For small genus, red usual first order Weil bound is the best one. Then from genus $g_2=\frac{\sqrt{q}(\sqrt{q}-1)}{2}$, yellow Ihara's second order Weil bound becomes the best, up to the genus $g_3=\frac{\sqrt{q}(q-1)}{\sqrt{2}}$ where green third order weil bound becomes the best. From some genus $g_4$, light blue fourth order Weil bound is the best up to some genus $g_5$, where dark blue fifth order Weil bound becomes better, and so on. Here, we have chosen~$q=2$,
for which the genera~$g_2,g_3,g_4,g_5$ are particularly small,
respectively about~$0.3, 1, 2.35$ and~$4.67$. This means for instance that the best bound for $q=2$ and $g=3, 4$ is the fourth order one!
\end{commentaire} \label{Graphes_N1} \label{Graphes_bornes_N1}
In order to illustrate the efficiency of our
algorithm, we display in section~\ref{table} a numerical
table.

We acknowledge that for any pair $(g, q)$ we have tested, a comparison  of our
numerical results with the table given on the website {\tt http://www.manypoints.org/} reveals that we
{\em always} recover the very same numerical upper bound for $N_q(g)$ than those coming from Oesterl\'e bounds! We were not been able to understand
this experimental observation.

Nevertheless, we think that the viewpoint introduced in this article is preferable to Osterl\'e's one. 
First, our viewpoint is more conceptual in nature. Any constraints we use to obtain our bounds come in a quite pleasant way either from algebraic-geometry or from arithmetic, as explained in the introduction of section~\ref{s_Absolute_bounds} in which we outline our approach. We think moreover that the graph displayed just above is very satisfying. 
Second, the viewpoint introduced in this article is perfectly adapted to the study of bounds for the numbers of rational points. As the reader can see, many\footnote{Other known results can be proved from this viewpoint, for instance the rationality of the Zeta function, using the vanishing of Hankel determinants. In order to save place, we have chosen to skip this.} known results can be understood with this viewpoint -even asymptotical ones, and new results can be proved. We are pretty sure that we have not extracted all potential outcomes of this viewpoint. 
Third, new questions can be raised from this viewpoint. We propose a few in Section~\ref{s_questions}.

To conclude Section~\ref{s_Absolute_bounds}, we push to the {\em infinite-order Weil bound} using the non-negativity of the Gram determinants of any orders, which imply that some symmetric matrix is semi-definite positive. Applied to some very simple vector, this lead us to Theorem~\ref{TVameliore}, a stronger form of \cite{T92} Tsfasman bound in that it gives a new interpretation for the {\em defect} of an exact tower as a limit of nice euclidean vectors in~$\left(\HVperp, \scalaire{\cdot}{\cdot}\right)$.

\medskip

In Section~\ref{s_Relative_bounds}, we study the relative situation. Given a finite morphism~$f : X \to Y$
between two absolutely irreducible smooth projective curves defined over~${\mathbb F}_q$, the pull-back functor on divisors from the bottom algebraic surface $Y\times Y$ to the top one $X \times X$ induces a map from the Neron Severi space $\NS(Y\times Y)_{\mathbb R}$ to $\NS(X\times X)_{\mathbb R}$.
Restricting this map to the subspaces~$\Frob_{X}$ and~$\Frob_{Y}$ generated by the classes of the graphs of iterations of the Frobenius morphisms,
and after a suitable normalization (which differs from the normalization
chosen in Section~\ref{s_Absolute_bounds}), we prove that it becomes an {\em isometric embedding}. We thus have an  orthogonal decomposition 
$\Frob_{X} = \Frob_{Y}^* \oplus \Frob_{X/Y}$,
involving a {\em relative subspace} $\Frob_{X/Y}$ of $\Frob_{X}$.

Now, Cauchy-Schwartz inequality applied to the orthogonal projections of $\Gamma^0_{X}$ and $\Gamma^1_{X}$ on the relative space~$\Frob_{X/Y}$ is equivalent to
the well known relative Weil  bound
 that $\vert \sharp X({\mathbb F}_q)
-  \sharp Y({\mathbb F}_q) \vert \leq 2(g_X-g_Y) \sqrt{q}$.

With regard to higher orders relative Weil  bounds, we encounter a difficulty
since the arithmetical constraints $\sharp X({\mathbb F}_{q^r}) -  \sharp Y({\mathbb F}_{q^r}) \geq \sharp X({\mathbb F}_q) -  \sharp Y({\mathbb F}_q)$ for~$r \geq 2$, similar to that used in Section\footnote{See the introduction of Section~\ref{s_Absolute_bounds}}~\ref{s_Absolute_bounds} does not hold true! Hence, the only constraints we are able to
use are the non-negativity of Gram determinants.
This lead to an inequality relating the
quantities $\sharp X({\mathbb F}_{q^r}) -  \sharp Y({\mathbb F}_{q^r})$
(corollary~\ref{RelatifNonTrivial}).
We also prove a bound involving four curves in a cartesian diagram under some smoothness assumption (Theorem~\ref{final}).

\medskip

As stated above, we end this article by a very short Section~\ref{s_questions}, in which we raise a few questions.

\section{The euclidean space~$\left(\HVperp, \scalaire{\cdot}{\cdot}\right)$}
\label{s_Euclidean_space}

Let $\Cun$ be an absolutely irreducible smooth projective curve defined over
the finite field $\FF_q$ with $q$ elements. We consider the Neron-Severi
space~$\NS(\Cun\times \Cun)_\RR = \NS(\Cun\times \Cun) \otimes_\ZZ \RR$ of the smooth 
surface~$\Cun\times \Cun$.
It is well known that the intersection product of two divisors~$\Gamma$
and~$\Gamma'$ on~$\Cun \times \Cun$ induces a symmetric bilinear
form on $\NS(\Cun\times \Cun)_\RR$, whose signature is given by
the following Hodge Index Theorem (e.g. \cite{Hartshorne}).

\begin{theo}[Hodge index Theorem]
The intersection pairing  is non-degenerate
on the space~$\NS(\Cun\times \Cun)_\RR$, and is definite negative on the
orthogonal supplement of any ample divisor.
\end{theo}

Let~$H = \Cun \times \{*\}$
and~$V = \{*\} \times \Cun$ be the horizontal and vertical classes
on~$\Cun \times \Cun$. Their intersection products are given
by~$H\cdot H = V\cdot V = 0$ and~$H\cdot V = 1$. The restriction
of the intersection product to~$\Vect(H,V)$ has thus signature~$(1,-1)$.
Moreover~$\Vect(H,V) \cap \Vect(H,V)^\perp = \{0\}$. By non-degeneracy, this gives rise to
an orthogonal
decomposition~$\NS(\Cun\times \Cun)_\RR
=
\HV \oplus \Vect(H,V)^\perp
$.
Let~$p$ be the orthogonal projection onto the non-trivial part~$\Vect(H,V)^\perp$, given by
\begin{equation}
\begin{array}{rccl}\label{projection}
p : & \NS(\Cun\times \Cun)_\RR &\longrightarrow &\Vect(H,V)^\perp \\
    & \Gamma                  &\longmapsto     & \Gamma - (\Gamma \cdot V)H - (\Gamma \cdot H)V.
\end{array}
\end{equation}
Since $H+V$ is ample,
the intersection pairing is definite negative on~$\Vect(H,V)^\perp$ by Hodge index Theorem.
Hence,~$\Vect(H,V)^\perp$ can be turned into an
euclidean space by defining a scalar product as the opposite of the intersection pairing.

\begin{defi} \label{def_scalar_product}
Let~$\HVperp = \Vect(H,V)^\perp$. 
We define on~$\HVperp$ a scalar product,
denoted by~$\scalaire{\cdot}{\cdot}$, as
$$
\scalaire{\gamma}{\gamma'}
= - \gamma \cdot \gamma',
\qquad \forall \gamma, \gamma' \in \HVperp.
$$
The associated norm on~$\HVperp$ is denoted by~$\|\cdot\|$.
\end{defi}

In the sequel of this article, all the computations will take place in the euclidean
space~$\left(\HVperp,\scalaire{\cdot}{\cdot}\right)$.
To begin with, 
let us compute this pairing between the iterates of the Frobenius morphism.

\begin{lem} \label{Intersections}
Let~$F : X \to X$ denotes the $q$-Frobenius morphism on $\Cun$, and~$F^k$ 
denotes the $k$-th iterate of~$F$ for any~$k\geq 0$, with the usual convention that~$F^0 = \Id_{\Cun}$.
We denote by~$\Gamma^{k}$ the Neron Severi class of the graph of~$F^k$. 
Then
\begin{equation*}\label{Gamma_k_Gamma_i}
\begin{cases}
\scalaire{p(\Gamma^{k})}{p(\Gamma^{k})}
=
2g q^k & \forall k\geq 0\\
\scalaire{p(\Gamma^{k})}{p(\Gamma^{{k+i}})}
= 
q^k\left((q^i+1) - \sharp X(\FF_q^i)\right)
& \forall k\geq 0, \, \forall i \geq 1.
\end{cases}
\end{equation*}
\end{lem}

\begin{preuve}
Since the morphism~$F^k$ is a regular map of degree~$q^k$, one
has~$\Gamma^k \cdot H=q^k$ and~$\Gamma^k \cdot V=1$.
Now, let $k \in \NN$ and $i \in \NN^*$. We consider the map
$\pi = F^k \times \Id : \Cun \times \Cun \rightarrow \Cun \times \Cun$, which sends $(P, Q)$ to $(F^k(P), Q)$. We have 
$\pi^*(\Delta)= \Gamma^{k}$ and $\pi_*\left(\Gamma^{{k+i}} \right) = q^k \Gamma^{i}$, so that by projection formula for the proper morphism $\pi$:
\begin{align*}
\Gamma^{{k}} \cdot \Gamma^{{k+i}} 
& = 
\pi^*(\Delta) \cdot \Gamma^{{k+i}}\\
&= 
\Delta \cdot \pi_*\left(\Gamma^{{k+i}} \right)\\ 
&= q^k \Delta \cdot \Gamma^{i}.
\end{align*}
If $i \geq 1$, then $\Delta \cdot \Gamma^{i} = \sharp \Cun(\FF_{q^i})$. If $i=0$, then $\Gamma^{0} = \Delta$ and
$\Delta \cdot \Gamma^{0} = \Delta^2=-2g$. The results follow since~$p(\Gamma^k) = \Gamma^k- H - q^k V$ by \eqref{projection}.
\end{preuve}

\begin{rk}
Of course, $\dim \HVperp \leq \dim NS(X\times X)_{\RR}-2 \leq 4g^2$. But as the reader can check, we are working along this article only on the subspace $\Frob$ of $\HVperp$ generated by the family $p(\Gamma^k), k\geq 0$. By (\cite{Zariski} chapter VII, appendix of Mumford), any trivial linear combination of this family is equivalent to the triviality of the same 
linear combination of the family $\varphi_{\ell}^k, k \geq 0$, of the iteration of the Frobenius endomorphism $\varphi_{\ell}$ on the Tate module $V_{\ell}(\hbox{Jac } X)$ for any prime $\ell \wedge q=1$. We deduce that for a given curve $X$ over $\FF_q$ of genus $g$, we actually are working in an euclidean space $\Frob_{X}$ of dimension equal to the degree of the minimal polynomial of $\varphi_{\ell}$ on $V_{\ell}(\hbox{Jac }  X)$.
\end{rk}


\section{Absolute bounds}\label{s_Absolute_bounds}
We observe, as a consequence of Lemma~\ref{Intersections}, that the non-negativity 
\begin{align*}
0 & \leq \Gram(p(\Gamma^{0}), p(\Gamma^{1})) =
\begin{vmatrix}
2g      &q+1-\sharp X(\FF_q)\\
q+1-\sharp X(\FF_q)    & 2gq
\end{vmatrix}\\
&=(2g\sqrt{q})^2-(q+1-\sharp X(\FF_q))^2
\end{align*}
is nothing else than Weil inequality. In this Section, we give other bounds using larger Gram determinants. For this purpose, we need preliminary notations, normalizations and results.

Normalizations in Section~\ref{SubsecGramMatrix} play two parts. They ease many formulae and calculations. They also make obvious that several features of the problem, such as the Gram matrix $\Gram(p(\Gamma^i) ; 0\leq i \leq n)$, or the forthcoming $n$-th Weil domains $\Dcal_n$, are essentially independent of $q$. The authors believe that these features deserve to be emphasized.

With regard to results,  let $n \in \NN^*$. To obtain the Weil bound of order $n$, let $X$ be an algebraic curve of genus $g$ defined over $\FF_q$. We use both geometric and arithmetic facets of $X$ as follows. We define from $X$ in~\eqref{Pn(X)} a point $P_n(X) \in \RR^n$, whose abscissa $x_1$ given in~\eqref{def_x_i} is essentially the {\em opposite} of $\sharp X(\FF_q)$, hence have to be {\em lower bounded}.
The algebraic-geometric facet of $X$ implies, thanks to Hodge index Theorem, that some Gram determinant involving the coordinates of $P_n(X)$ is non-negative. This is traduced on the point $P_n(X)$ in Lemma~\ref{T(X)DansDn} in that it do lies inside some convex {\em $n$-th Weil domain} $\overline{\Dcal}_n$, studied in Subsection~\ref{SubsecDn}. The arithmetic facet of $X$ is used via {\em Ihara constraints}, that for any $r \geq 2$ we have
$\sharp X(\FF_{q^r}) \geq \sharp X(\FF_q)$. We then define in Subsection~\ref{SubsecCurveLocus}  the convex {\em $n$-th Ihara domain $\Hcal_n^g$ in genus $g$}. We state in Proposition~\ref{T(X)DansDncapHng} that for any curve $X$ of genus $g$ defined over $\FF_q$, we have $P_n(X) \in  \overline{\Dcal}_n \cap \Hcal_n^g$. We are thus reduced to minimize the convex function $x_1$ on the convex domain $ \overline{\Dcal}_n \cap \Hcal_n^g$. To this end, we establish in Subsection~\ref{SubsecCurveLocus}
the optimization criteria~\ref{criteria_min_x_1}, decisive for the derivations of higher order Weil bounds.

This criteria is used to prove new bounds in Subsection~\ref{NewBounds}. Finally, we obtain in Subsection~\ref{asymptotique} a refined form of Tsfasman upper bound for asymptoticaly exact towers.



\subsection{The Gram matrix of the normalized classes of iterated Frobenius} \label{SubsecGramMatrix}

In this Subsection, we choose to normalize the Neron-Severi classes of the iterated Frobenius morphisms
as follows.
For any $n \geq 1$, and any curve $X$ defined over $\FF_q$ of genus $g \neq 0$, we define a point $P_n(X) \in \RR^n$ whose $x_i$-coordinates is very closely related to $\sharp X(\FF_{q^i})$.

\begin{defi} \label{Normalization}
Let $X$ be an absolutely irreducible smooth projective curve defined over the finite field $\FF_q$. For~$k \in \NN$, we put
\begin{equation}\label{label gammak}
\gamma^k = \frac{1}{\sqrt{2g q^k}} p(\Gamma^k) \in \HVperp,
\end{equation}
so that by Lemma~\ref{Intersections} we have~$\scalaire{\gamma^{k}}{\gamma^{k}} = 1$. For~$i\in\NN^*$,
we put
\begin{equation}\label{def_x_i}
x_i
\overset{\text{def.}}{=}
\scalaire{\gamma^{k}}{\gamma^{{k+i}}}
\overset{\text{lemma~\ref{Intersections}}}{=} \frac{(q^i+1) - \sharp\Cun(\FF_{q^i})}{2g\sqrt{q}^i}.
\end{equation}
%
%
%
For any $n \geq 1$, we define the point
\begin{equation}Ê\label{Pn(X)}
P_n(X) = (x_1, \ldots, x_n) \in \RR^n.
\end{equation}
\end{defi}

\begin{rk}
Note first that $X$ is Weil-maximal if and only if $x_1=-1$, and is Weil minimal if and only if $x_1=1$. Second, that to give {\em upper} bounds for $\sharp X(\FF_q)$ amount to give {\em lower} bounds for $x_1$.
\end{rk}

From this Definition~\ref{Normalization} and Lemma~\ref{Intersections}, the Gram matrix of the family~$\gamma^0,\ldots,\gamma^n$ in~$\HVperp$ is
\begin{equation}
 \Gram(\gamma^0,\ldots,\gamma^n)
=
\begin{pmatrix} \label{MatrGram}
1      & x_1    & \cdots & x_{n-1} & x_n\\
x_1    & \ddots & \ddots &        & x_{n-1} \\
\vdots & \ddots & \ddots & \ddots & \vdots \\
x_{n-1} &        & \ddots & \ddots & x_1\\
x_n    & x_{n-1} & \cdots & x_1     & 1
\end{pmatrix}
\end{equation}

\subsection{Some identities involving Toeplitz matrices}
The main result of this Subsection is Lemma~\ref{lem_factorization}, providing the existence of the factorization in Definition~\ref{def_minors}, in which the $G_n^-$ factor plays a fundamental part in the following of the article.

\bigskip

A {\em symmetric Toeplitz matrix} is a symmetric matrix
whose entries $x_{i, j}$ depend only on $\vert i-j \vert$, that is are constant along the diagonals parallel to the main
diagonal \cite[\S0.9.7]{MatrixAnalysis}, that is
\begin{equation} \label{Tn+1}
T_{n+1}(x_0,x_1,\ldots,x_n)
=
\begin{pmatrix}
x_0      & x_1    & \cdots & x_{n-1} & x_n\\
x_1    & \ddots & \ddots &        & x_{n-1} \\
\vdots & \ddots & \ddots & \ddots & \vdots \\
x_{n-1} &        & \ddots & \ddots & x_1\\
x_n    & x_{n-1} & \cdots & x_1     & x_0
\end{pmatrix}
\end{equation}
The symmetric Toeplitz matrix is said to be {\em normalized} if~$x_0 = 1$.

An {\em Hankel matrix} is a matrix
whose entries $x_{i, j}$ depend only on $i+j$, that is are constant along the anti-diagonals parallel to the main
anti-diagonal \cite[\S0.9.8]{MatrixAnalysis}, that is
\begin{equation} \label{Hn+1}
H_{n+1}(x_0,x_1,\ldots,x_{2n})
=
\begin{pmatrix}
x_0    & x_1    & x_2    & \cdots  & x_{n}\\
x_1    & \reflectbox{$\ddots$} &  & \reflectbox{$\ddots$}  & x_{n+1} \\
x_2    &  & \reflectbox{$\ddots$} &\reflectbox{$\ddots$}  & \vdots \\
\vdots & \reflectbox{$\ddots$}  & \reflectbox{$\ddots$} &  & x_{2n-1}\\
x_{n} & x_{n+1} & \cdots & x_{2n-1} & x_{2n}
\end{pmatrix}
\end{equation} 

\begin{rk}
It should be noticed that in this article, any matrix is indexed by its size, while this doesn't hold for determinants, for instance see  formulae~\eqref{G2} and ~\eqref{G3}.
\end{rk}

\begin{lem} \label{lem_factorization}
For any~$n\geq 1$, we abbreviate~$T_n(x_0,\ldots,x_{n-1})$ by~$T_{n}$.
The determinant of this Toeplitz matrix factorizes as
\begin{equation*}
\Det\left(T_{2n}\right)
=
\Det\left(T_n + H_{n}(x_{2n-1},\ldots,x_1)\right)
\times
\Det\left(T_{n} - H_{n}(x_{2n-1},\ldots,x_1)\right)
\end{equation*}
and
\begin{equation*}
\Det\left(T_{2n+1}\right)
=
\Det
\begin{pmatrix}
T_{n} + H_{n}(x_{2n},\ldots,x_2) & {}^t X_n \\
2X_n & x_0
\end{pmatrix}
\times
\Det\left(T_{n} - H_{n}(x_{2n},\ldots,x_2)\right),
\end{equation*}
where~$X_n = (x_n,\ldots,x_1)$.
\end{lem}

\begin{preuve}
{\em Even size}. The operations on the columns~$C_j \leftarrow C_j + C_{2n-(j-1)}$
for~$1\leq j\leq n$, followed by the operations on the
lines~$L_{n + i} \leftarrow L_{n + i} - L_{n-(i-1)}$ for~$1\leq i\leq n$ lead
to the determinant of a $2\times 2$ blocks upper triangular matrix
whose value is the expected product.

{\em Odd size.} The operations on the columns~$C_j \leftarrow C_j + C_{2n+1-(j-1)}$
for~$1\leq j\leq n$, followed by the operations on the
lines~$L_{n+1 + i} \leftarrow L_{n+1 + i} - L_{n+1-i}$ for~$1\leq i\leq n$ lead
to the result in the same manner.
%
\end{preuve}

\begin{defi}\label{def_minors}
For~$n=0$, we put~$G_0=G_0^-=G_0^+=1$, and for~$n\geq 1$, we put
$$
G_n(x_1,\ldots,x_n)
\overset{\text{def.}}{=}
\Det(T_{n+1}(1,x_1,\ldots,x_n))
\overset{\text{Lemma~\ref{lem_factorization}}}{=}
G_n^-(x_1,\ldots,x_n) \times G_n^+(x_1,\ldots,x_n),
$$
where the last factorization is the one proved in Lemma~\ref{lem_factorization}
in the same order.
%
\end{defi}

\begin{rk}
The determinant~$G_n$ is the determinant
of a matrix of size~$n+1$. Note that for $n=2m +1$ odd, both $G_n^-$ and $G_n^+$ are polynomials of degree $m+1$ in $x_1, \ldots,  x_n$, while for $n=2m$ even, $G_n^-$ has degree $m+1$, while $G_n^+$ have degree $m$. 
We will see later in the article that Weil bound of order n for $\sharp X(\FF_q)$ depends heavily on the hypersurface $\{G_n^-=0\}$.
\end{rk}

In order to raise any ambiguity and for later purpose,
let us write down these  determinants for~$n=1,2$ and~$3$:
\begin{align}
G_1(x_1)
&=
\begin{vmatrix}
1&x_1\\
x_1&1
\end{vmatrix}
=
\underbrace{(1+x_1)}_{G_1^-(x_1)}\times
\underbrace{(1-x_1)}_{G_1^+(x_1)}\label{G1}\\
G_2(x_1,x_2)
&=
\begin{vmatrix}
1&x_1&x_2\\
x_1&1&x_1\\
x_2&x_1&1
\end{vmatrix}
=
\underbrace{
\begin{vmatrix}
1+x_2&x_1\\
x_1+x_1&1
\end{vmatrix}
}_{G_2^-(x_1,x_2)}
\times
\underbrace{(1-x_2)}_{G_2^+(x_1,x_2)}
\label{G2}\\
G_3(x_1,x_2,x_3)
&=
\begin{vmatrix}
1&x_1&x_2&x_3\\
x_1&1&x_1&x_2\\
x_2&x_1&1&x_1\\
x_3&x_2&x_1&1
\end{vmatrix}
=
\underbrace{
\begin{vmatrix}
1+x_3&x_1+x_2\\
x_1+x_2&1+x_1
\end{vmatrix}
}_{G_3^-(x_1,x_2,x_3)}
\times
\underbrace{
\begin{vmatrix}
1-x_3&x_1-x_2\\
x_1-x_2&1-x_1
\end{vmatrix}
}_{G_3^+(x_1,x_2,x_3)}\label{G3}
\end{align}

\begin{lem} \label{lem_2G_n}
For any~$n\geq 1$ and~$\epsilon = \pm$ or nothing, we
abbreviate~$G_n^\epsilon(x_1,\ldots,x_n)$ by~$G_n^\epsilon$. Then one has
$$
2G_n
=
G_{n-1}^+ \times G_{n+1}^-
+
G_{n-1}^-\times G_{n+1}^+.
$$
\end{lem}

\begin{preuve}
Let~$S_{n+1}^-$ be the matrix obtained from~$T_n$
by adding~${}^t(x_n,\ldots,x_1)$ to the first column,
and~$S_{n+1}^+$ be the matrix obtained from~$T_n$
by removing~${}^t(x_n,\ldots,x_1)$ to the first column.
Then by multilinearity, there exists a polynomial~$R_n$ such that:
$$
\Det\left(S_{n+1}^-\right) = G_n + R_n
\quad\text{and}\quad
\Det\left(S_{n+1}^+\right) = G_n - R_n.
$$
On the other hand, using the same transformations as
in the proof of Lemma~\ref{lem_factorization}, one get
$$
\Det(S_{n+1}^-) = G_{n-1}^+ G_{n+1}^-
\quad\text{and}\quad
\Det(S_{n+1}^+) = G_{n-1}^- G_{n+1}^+
$$
Adding $\Det(S_{n+1}^-)$ and $\Det(S_{n+1}^+)$  allows us to conclude.
\end{preuve}

\subsection{The domain of positive definite normalized Toeplitz matrices} \label{SubsecDn}
In this Subsection, we study for any $n \geq 1$ some domain $\Dcal_n \subset \RR^n$,
whose closure $\overline{\Dcal}_n$ also plays a fundamental part in this article.

\bigskip

\begin{defi}
Let $n \in {\mathbb N}^*$. We denote by $\Dcal_n$ the set of $(x_1, \ldots, x_n) \in \RR^n$, such that
the symmetric normalized Toeplitz matrix $T_{n+1}(1, x_1, \ldots, x_n)$ is positive definite.
\end{defi}


The domain~$\Dcal_n$ can be characterized in several useful ways.

\begin{prop} \label{DnRecursif}
For~$n \geq 1$, the domain~$\Dcal_n$ is a convex subset of~$\left]-1,1\right[^n$, which can be written as follows.
\begin{enumerate}
\item \label{eq_def_D_n} First, one has:
$$
\Dcal_n
=
\left\{
(x_1,\ldots,x_n) \in \RR^n \mid G_i(x_1,\ldots,x_i) > 0,\; \forall \; 1\leq i \leq n
\right\}.
$$
\item \label{def_Dn_rec} Recursively,~$\Dcal_1 = \left]-1,1\right[$ and for any~$n\geq 2$
\begin{align*}
\Dcal_n
&=
\left\{
(x_1,\ldots,x_{n-1},x_n) \in \Dcal_{n-1} \times \RR  \mid
G_n(x_1,\ldots,x_n) > 0
\right\} \\
&=
\left\{
(x_1,\ldots,x_{n-1},x_n) \in \Dcal_{n-1} \times \RR \mid
G_n^-(x_1,\ldots,x_n) > 0
\text{ and } G_n^+(x_1,\ldots,x_n) > 0
\right\}.
\end{align*}
\item \label{def_Dn_graphe} $\Dcal_n$ is also the set of points between the graphs of two functions
from~$\Dcal_{n-1}$ to~$\RR$: for~$\varepsilon = \pm$, there exists
a polynomial~$\widetilde{G}_n^\varepsilon \in \QQ[x_1,\ldots,x_{n-1}]$, such
that~$G_n^\epsilon = -\epsilon G_{n-2}^\epsilon x_n
+ \epsilon \widetilde{G}_n$ and
\begin{align*}
\Dcal_n
=
\left\{
(x_1, \ldots, x_{n-1}, x_n) \in \Dcal_{n-1}\times \RR \mid
\frac{\widetilde{G}^-_n(x_1,\ldots,x_{n-1})}{G_{n-2}^-(x_1,\ldots,x_{n-2})}
< x_n <
\frac{\widetilde{G}_n^+(x_1,\ldots,x_{n-1})}{G_{n-2}^+(x_1,\ldots,x_{n-2})}
\right\}.
\end{align*}
\end{enumerate}
\end{prop}

\begin{preuve} To prove the convexity, we remark that
both sets of normalized symmetric Toeplitz matrices and of symmetric positive definite matrices are convex,
so that $\Dcal_n$ is convex.
Moreover,~$\Dcal_n \subset \left]-1,1\right[^n$
because if~$T_{n+1}(1,x_1,\ldots,x_n)$ is positive definite, then all its
principal minors are positive. In
particular, for any $1 \leq i \leq n$, the $2\times 2$
minors~$\begin{vsmallmatrix}1&x_i\\x_i&1\end{vsmallmatrix}$  is positive.

Item~\eqref{eq_def_D_n}  follows from the well known fact that an $n\times n$ symmetric matrix is definite positive
if and only if all its $n$
leading principal minors (obtained by deleting the $i$ last rows 
and columns for~$0\leq i\leq n-1$) are positive 
\cite[Theorem 7.2.5]{MatrixAnalysis}.
Hence, the first characterization of item~\eqref{def_Dn_rec} follows from item~\ref{eq_def_D_n} . To prove the second characterization of item~\eqref{def_Dn_rec},
we use Lemma~\ref{lem_2G_n} stating
that~$2G_{n-1} = G_{n-2}^-G_n^+ + G_{n-2}^+G_n^-$. By induction, this allow
us to prove that both factors are positive if, and only if, 
the product
$$
G_n(x_1,\ldots,x_n)=G_n^-(x_1,\ldots,x_n)G_n^+(x_1,\ldots,x_n)
$$
is positive for any~$(x_1,\ldots,x_{n-1}) \in \Dcal_{n-1}$.

To prove item~\eqref{def_Dn_graphe}, for~$\epsilon = \pm$, formulae
defining the
polynomials~$G_n^\epsilon$ in Lemma~\ref{lem_factorization} and
Definition~\ref{def_minors} imply, developing along their first column and taking advantage of the very particular forms~\eqref{Tn+1} and~\eqref{Hn+1}, that both
polynomials
have degree~$1$ in~$x_n$, of the form
$$G_n^\epsilon = -\epsilon G_{n-2}^\epsilon x_n + \varepsilon \widetilde{G}_n^\epsilon(x_1,\ldots,x_{n-1})$$
 for some~$\widetilde{G}_n^\epsilon \in \QQ[x_1,\ldots,x_{n-1}]$. Since by
item~\ref{def_Dn_rec} we have $G_{n-2}^{\varepsilon}>0$, the set~$\Dcal_n$ is thus equal to the set
of~$(x_1,\ldots,x_{n-1},x_n)\in\Dcal_{n-1} \times \RR$, such that:
$$
\frac{\widetilde{G}^-_n(x_1,\ldots,x_{n-1})}{G_{n-2}^-(x_1,\ldots,x_{n-2})}
< x_n <
\frac{\widetilde{G}_n^+(x_1,\ldots,x_{n-1})}{G_{n-2}^+(x_1,\ldots,x_{n-2})},
$$
and the proof is complete.
\end{preuve}

The following Proposition is useful for later purpose in Section~\ref{NewBounds}, where convexity plays a quite important part.

\begin{prop} \label{ConvexeConcave}
The locus~$\{(x_1,\ldots,x_n) \in \Dcal_{n-1}\times\RR \mid G_n^-(x_1,\ldots,x_n) >0\}$ is convex,
while the locus~$\{(x_1,\ldots,x_n) \in \Dcal_{n-1}\times\RR \mid G_n^+(x_1,\ldots,x_n) >0\}$ is concave.
\end{prop}

\begin{preuve} By item~\ref{def_Dn_graphe} of Proposition~\ref{DnRecursif},
$\Dcal_n$ is the set
of~$(x_1,\ldots,x_{n-1},x_n) \in \Dcal_{n-1}\times\RR$ such that~$x_n$ is 
between two functions from~$\Dcal_{n-1}$ to~$\RR$. By convexity of $\Dcal_n$ from Proposition~\ref{DnRecursif}, the
lower function must be concave has a function on $\Dcal_{n-1}$ while the upper one must be convex. The Proposition
follows easily.
\end{preuve}

The closure of~$\Dcal_n$ is given by the following Proposition.

\begin{prop} \label{prop_carac_Dn_bar}
The closure~$\overline{\Dcal}_n$ of~$\Dcal_n$ in~$\RR^n$ corresponds to the
set of $(n+1)\times(n+1)$ normalized symmetric Toeplitz matrices
which are positive semi definite. Moreover, we have
$$
\overline{\Dcal}_n
=
\left\{
(x_1,\ldots,x_n) \in \RR^n \mid G_{n,I}(x_1,\ldots,x_i) \geq 0,\; \forall I \subset \{0, 1,\ldots,n\}
\right\},
$$
where, for $I \subset \{0, 1,\ldots,n\}$, we denote by
$G_{n,I}(x_1,\ldots,x_n)$ the principal minor
of the normalized symmetric Toeplitz
matrix~$T_{n+1}(1,x_1,\ldots,x_n)$ obtained by deleting the lines and
columns whose indices are not in~$I$. 
\end{prop}

\begin{preuve}
This is a consequence of the fact that a matrix
is {\em positive semi definite} if and only if all the {\em principal minors} (the ones obtained by deleting the same
subset of lines and columns) are non-negative (i.e.~$\geq 0$) \cite[last exercise after Theorem 7.2.5]{MatrixAnalysis}.
%
%
\end{preuve}

\subsection{The curves locus} \label{SubsecCurveLocus}
In this Subsection, we introduce Ihara's constraints and the resulting {\em $n$-Ihara domain $\Hcal_n^g$ for genus g}. We then gather the key results for the derivation of higher order Weil bounds. The first one is Proposition~\ref{T(X)DansDncapHng}, the second one is criteria~\ref{criteria_min_x_1}. 

\bigskip


\subsubsection{The $n$-the Weil domain $\overline{\Dcal}_n$}

\begin{lem} \label{T(X)DansDn}  
Let $X$ be an absolutely irreducible smooth projective curve defined
over~$\FF_q$. Then the point~$P_n(X)= (x_1, \ldots, x_n)$,
 as defined in~\eqref{Pn(X)}, belongs to~$\overline{\Dcal}_n$.
\end{lem}

\begin{preuve}
By~\eqref{def_x_i}, we have
$x_i= \scalaire{\gamma^{k}}{\gamma^{k+i}}$ where the $\gamma^k$ are defined in~\eqref{label gammak}, so that it is easily seen that for any $I \subset \{0, 1,\ldots,n\}$, we have
$$
G_{n,I}(x_1,\ldots,x_n)= \Gram(\gamma^i, \, i \notin I),
$$
for $G_{n,I}$ s defined in Proposition~\ref{prop_carac_Dn_bar}. The non-negativity of $G_{n, I}$ follows as a Gram determinant in an euclidean space, hence the Lemma by Proposition~\ref{prop_carac_Dn_bar}.
\end{preuve}

\begin{defi}
$\overline{\Dcal}_n$ is called the $n$-th Weil domain.
\end{defi}

As a first illustration of the informations contained in these Weil domains, we prove the following simple result. From~\eqref{G2}, the non-negativity of the Gram determinants~$G_2^-$
writes~$x_2\geq 2x_1^2-1$. Taking~\eqref{def_x_i} into account, this gives immediately

\begin{prop} \label{X(Fq2)}
For any curve $X$ of genus $g \neq 0$, we have
\begin{equation*}
\sharp X({\mathbb F}_{q^2}) -(q^2+1) \leq 2gq-\frac{1}{g}\Bigl(\sharp X({\mathbb F}_q)-(q+1)\Bigr)^2.
\end{equation*}
\end{prop}

This Proposition means that, for a given non rational curve $X$, any lower bound for the deviation
of $\sharp X({\mathbb F}_q)$ to $q+1$ yields to a better upper bound than Weil's one for $\sharp X({\mathbb F}_{q^2})$. In the same way, for any given order $n$, the non-negativity of the $3\times 3$ determinant $\Gram(\gamma^0, \gamma^1, \gamma^n)$ gives a quite ugly upper bound of similar nature for $\sharp X({\mathbb F}_{q^n})$ in terms of $\sharp X({\mathbb F}_{q})$ and $\sharp X({\mathbb F}_{q^{n-1}})$.

\begin{rk}
Note that this Proposition is a refinement of the following well known particular case: if $X$ is either Weil maximal or minimal over ${\mathbb F}_q$, then $\sharp X({\mathbb F}_q)-(q+1)=\pm2g\sqrt{q}$, and this Proposition asserts that then $X({\mathbb F}_{q^2}) -(q^2+1) \leq 2gq-\frac{4g^2q}{g}=-2gq$, so that $X$ is Weil minimal over ${\mathbb F}_{q^2}$. Note also that curves such that this inequality is an equality are those curves such that the corresponding point $P_2(X)=(x_1, x_2)\in \overline{\Dcal}_2$ lies on the bottom parabola $x_2=2x_1^2-1$ of the Ihara domain drawn in figure~\ref{figure_Ihara_domain} below. The above particular case is that of curves corresponding to the corner points $(-1, 1)$ and $(1, 1)$.
\end{rk}

\subsubsection{Ihara constraints : the $n$-th Ihara domain $\Hcal_n^g$ in genus $g$}

If~$(x_1,\ldots,x_n) \in\RR^n$ comes from a curve $X$ over $\FF_q$ by formulae~\eqref{def_x_i}, 
then we have seen in Lemma~\ref{T(X)DansDn} that
~$(x_1,\ldots,x_n)$ lies on  the closure~$\overline{\Dcal}_n$ 
of~$\Dcal_n$. But there are also other constraints, resulting
from the arithmetical inequalities $\sharp X(\FF_{q^i}) \geq \sharp X(\FF_q)$ for any $i\geq 1$. These
inequalities, by~(\ref{def_x_i}), write 
\begin{equation} \label{Inegalites_x_i}
\forall i\geq 2,\qquad
x_i \leq \frac{x_1}{q^{\frac{i-1}{2}}} + \frac{q^{i-1}-1}{2g q^{\frac{i-2}{2}}}.
\end{equation}
Let
\begin{equation} \label{alpha}
\alpha = \frac{1}{\sqrt{q}}.
\end{equation}
For~$g\geq 0$, and~$n\geq 2, i\geq 2$, we define
\begin{equation} \label{EquationsToits}
h_i^g(x_1,x_i)
=
x_i - \alpha^{i-1} x_1 - \frac{1}{2g\alpha}
\left(\frac{1}{\alpha^{i-1}} - \alpha^{i-1}\right).
\end{equation}
Then it is easily seen that (\ref{Inegalites_x_i}) is equivalent to
\begin{equation} \label{P(x_1)}
h_i^g(x_1,x_i)
\leq 0.
\end{equation}
We now define

\begin{defi} Let $n \geq 1$ and $g \geq 0$. We define:
\begin{enumerate}
\item The $n$-th Ihara domain in genus $g$ as
\begin{equation} \label{defHng}
\Hcal_n^g
=
\left\{(x_1,\ldots,x_n) \in \RR^n
\mid
h_i^g(x_1,x_i) \leq 0,\; \forall \; 2\leq i\leq n\right \} \
\end{equation}
\item The $n$-th Ihara line in genus $g$ as
\begin{equation} \label{defLng}
\Lcal_n^g
=
\left\{(x_1,\ldots,x_n) \in \RR^n
\mid
h_i^g(x_1,x_i) = 0,\; \forall \; 2\leq i\leq n\right \}. 
\end{equation}
\end{enumerate}
\end{defi}

Each~$\Lcal_n^g$ is a line which is identified with~$\RR$ using
the first coordinate~$x_1$ as a parameter. We denote by
\begin{equation} \label{P(x_1)}
P_n^g(x_1) = \Bigl(x_1, \alpha x_1 + \frac{q-1}{2g}, \alpha^2x_1+\frac{q^2-1}{2g\sqrt{q}}, \ldots, \alpha^{n-1}x_1+ 
\frac{q^{n-1}-1}{2g q^{\frac{n-2}{2}}} \Bigr) 
\end{equation}
the point of $\Lcal_n^g$ with parameter $x_1$.

\medskip

In the same way, for~$g=\infty$, we define the {\em $n$-th Ihara infinite line} by
$$
\Lcal_n^\infty
=
\left\{
(x_1,\ldots,x_n) \in \RR^n
\mid
x_i = \alpha^{i-1}x_1,\; 2\leq i\leq n
\right\}.
$$
It follows from this the following Proposition.

\begin{prop} \label{T(X)DansDncapHng}
Let $X$ be an absolutely irreducible smooth projective curve defined over the finite field $\FF_q$ and $n\geq 2$. Then
$P_n(X)=(x_1, \ldots, x_n) \in  \overline{\Dcal}_n \cap \Hcal_n^g$, where the $x_i$'s are defined from $X$ by~\eqref{def_x_i}.
\end{prop}

\bigskip

From this Proposition~\ref{T(X)DansDncapHng} and the remark below Definition~\ref{def_minors}, the strategy is the following.

\begin{quote}
{\em{\bfseries Strategy.~---} For each~$n\geq 1$, we minimize the
coordinate function~$x_1(P)$ for $P$ lying inside the compact convex domain~$\overline{\Dcal}_n \cap \Hcal_n^g$. By Proposition~\ref{T(X)DansDncapHng}, this leads to a lower bound for $x_1(P_n(X))$ for any curve $X$ of genus $g$ over $\FF_q$, hence by~(\ref{def_x_i})
 to an upper bound for the number~$\sharp X(\FF_q)$. It turns that for given $q$ and genus $g$,
this bound is better and better for larger and larger $n$, up to an optimal one.}
\end{quote}

We are thus face to an optimization problem. Since the domain $\overline{\Dcal}_n \cap \Hcal_n^g$
on which we have to minimize the convex function~$x_1$ is also convex
 by Proposition~\ref{DnRecursif}, one has the following necessary and sufficient characterization
of the minimum from  \cite[Th\'eor\`eme~2.2]{Urruty} where in our case the {\em active constraints} are $G_n^-=0$ and
$h_i^g=0$ for $2 \leq i \leq n$, and where $h_i^g$ are defined in~\eqref{EquationsToits}. 

\medskip

{\it
Let~$P_0 \in \overline{\Dcal}_n \cap \Lcal_n^g$. Suppose that:
\begin{itemize}
\item for any $I \subset \{1,\ldots,n+1\}$ and $\varepsilon = \pm$ such that $G^{\varepsilon}_{n, I}(P_0)=0$,
there exists~$\mu^{\varepsilon}_I\geq 0$,
\item for any $2 \leq i \leq n$, there exists $\mu_i \geq 0$,
\end{itemize}
 such that
\begin{equation}\label{carac_min}
\nabla x_1 (P_0)
-
\sum_{I \subset \{1,\ldots,n+1\}} \mu^-_I \nabla G^-_{n,I}(P_0)
-
\sum_{I \subset \{1,\ldots,n+1\}} \mu^+_I \nabla G^+_{n,I}(P_0)
+
\sum_{i=2}^n \mu_i \nabla h_i^g(P_0) = 0.
\end{equation}
Then,~$x_1(P_0) = \min \{x_1 \mid(x_1,\ldots,x_n)\in\Dcal_n \cap \Hcal_n^g\}$.
}
\medskip

We deduce from that the following criteria.

\begin{criteria}[for minimizing~$x_1$]\label{criteria_min_x_1}
If~$P_0 = (x_1,\ldots,x_n) \in\overline{\Dcal}_n \cap \Lcal_n^g$ satisfies
\begin{enumerate}
\item $G_n^-(P_0) = 0$;
\item $\partial_i G_n^-(P_0) \geq 0$  for all $2 \leq i \leq n$;
\item $\sum_{i=1}^n \partial_i G_n^-(P_0) \alpha^{i-1} > 0$,
\end{enumerate}
then~$P_0$ minimizes the first coordinate~$x_1$
on~$\overline{\Dcal}_n \cap \Hcal_n^g$.
\end{criteria}

\begin{rk}
Corollary~\ref{DeriveeDirectionnelle} bellow states that the first requirement implies the third one. 
\end{rk}

\begin{preuve}
Suppose that the assumptions of the criteria hold true.  Then the system
$$
\begin{cases}
\left(\sum_{i=1}^n \partial_i G_n^-(P_0) \alpha^{i-1}\right)\mu^-_{\emptyset} = 1\\
\mu_2 = \partial_2 G_n^-(P_0) \mu^-_{\emptyset}\\
\text{\hspace{0.62cm}}\vdots\\
\mu_n = \partial_n G_n^-(P_0) \mu^-_{\emptyset}\\
\end{cases}
$$
has a solution $\mu^-_{\emptyset} > 0, \mu_2\geq 0, \ldots, \mu_n \geq 0$. As easily seen, this solution also
satisfy 
$$
\begin{pmatrix}
1\\0\\\vdots\\\vdots\\0
\end{pmatrix}
-
\mu^-_{\emptyset}
\begin{pmatrix}
\partial_1 G_n^-(P_0)\\
\partial_2 G_n^-(P_0)\\
\vdots\\\vdots\\
\partial_n G_n^-(P_0)
\end{pmatrix}
+
\mu_2
\begin{pmatrix}
-\alpha\\1\\0\\\vdots\\0
\end{pmatrix}
+
\cdots
+
\mu_n
\begin{pmatrix}
-\alpha^{n-1}\\0\\\vdots\\0\\1
\end{pmatrix}
=
\begin{pmatrix}
0\\\vdots\\\vdots\\0
\end{pmatrix},
$$
hence~(\ref{carac_min}) holds with the choice $\mu_I=0$ for any $I\neq \emptyset$ and $\mu^+_{\emptyset}=0$, 
so that~$x_1(P_0) = \min \{x_1 \mid(x_1,\ldots,x_n)\in \overline{\Dcal}_n \cap \Hcal_n^g\}$, and the criteria is
proved.
\end{preuve}

\subsection{Weil bound and Ihara bound as bounds of order $1$ and $2$} \label{KnownBounds}
We prove here that our viewpoint enables to deduce Weil bound and Ihara bound. We think that this Subsection
is of interest since it shows on simple cases how this method works, especialy in the case of Ihara bound thanks to figure~\ref{figure_Ihara_domain}.

\subsubsection{Weil bound as a bound of order $1$} \label{s_Weil_bound}

With this viewpoint, the usual Weil bound comes from Cauchy-Schwartz
inequality applied to~$\gamma^0$ and~$\gamma^1$. Indeed, one  have by~(\ref{G1})
$$
\Gram(\gamma^0,\gamma^1) = \begin{vmatrix}1 & x_1\\x_1&1\end{vmatrix} \geq 0,
\qquad\hbox{that is}\qquad
\vert x_1\vert \leq 1.
$$
By~(\ref{def_x_i}), we have just recovered the Weil bound
\begin{equation}\tag{\bfseries Weil ``first order" bound}\label{Weil_bound}
\boxed{
\vert \sharp \Cun(\FF_q) -(q+1)\vert \leq 2g\sqrt{q}
}
\end{equation}

Just for fun, we can also easily recover the fact that a curve which is maximal over~$\FF_q$
must be minimal over the~$\FF_{q^{2i}}$ and maximal of the~$\FF_{q^{2i+1}}$.
Indeed, being maximal over~$\FF_q$ means by~(\ref{def_x_i}) that~$x_1 = -1$, therefore
$$
\Gram(\gamma^0,\gamma^1,\gamma^2)
=
\begin{vmatrix}
1&-1&x_2\\
-1&1&-1\\
x_2&-1&1
\end{vmatrix}
=
-(1-x_2)^2 \geq 0,
\qquad\hbox{so that}\qquad
x_2 = 1,
$$
that is by~(\ref{def_x_i}) that~$X$ is minimal over~$\FF_{q^{2}}$. Then
$$
\Gram(\gamma^0,\gamma^2,\gamma^3)
=
\begin{vmatrix}
1&1&x_3\\
1&1&-1\\
x_3&-1&1
\end{vmatrix}
=
-(1+x_3)^2 \geq 0,
\qquad\hbox{so that}\qquad
x_3 = -1,
$$
that is~$X$ is maximal over~$\FF_{q^{3}}$, and so on\dots In the same way,
if~$X$ is minimal over~$\FF_q$, that is if~$x_1 = 1$, then
$$
\Gram(\gamma^0,\gamma^1,\gamma^2)
=
\begin{vmatrix}
1&1&x_2\\
1&1&1\\
x_2&1&1
\end{vmatrix}
=
-(x_2-1)^2 \geq 0,
\qquad\hbox{that is}\qquad
x_2 = 1,
$$
and~$X$ still minimal over~$\FF_{q^2}$. And so on again...

\subsubsection{Ihara  bound as a bound of order~$2$}\label{s_Ihara_bound}

\begin{figure}
\begin{center}
\begin{tikzpicture}[scale=1.5,>=latex]
\draw [very thin, gray!20,step = 1] (-1.2,-1.2) grid (1.2,1.2);
\draw[color=gray!20,draw=black,very thick] (-1,1) -- (1,1) -- plot [domain=1:-1] (\x,{2*(\x)^2-1}) -- cycle;
\filldraw[color=gray!20,draw=black,very thick]
(1,1) -- (1/2,1) -- plot[domain = 1/2 : {(1-57^(1/2))/8}] (\x, 0.5*\x + 0.75) --  plot [domain={(1-57^(1/2))/8}:1] (\x,{2*(\x)^2-1}) -- cycle;
\draw [very thick] (-1.2,{-1.2/2+3/4}) -- (1.2,{1.2/2+3/4}) node[below, anchor=west]
{$\Lcal_2^g, \; g> \frac{\sqrt{q}(\sqrt{q}-1)}{2}$} ;
\draw [very thick, gray] (-1.2,{-1.2/2+3/2}) -- (1.2,{1.2/2+3/2}) node[black,below, anchor=west]
{$\Lcal_2^g, \; g = \frac{\sqrt{q}(\sqrt{q}-1)}{2}$} ;
\draw [very thick, gray!40] (-1.2,{-1.2/2+1.75}) -- (1.2,{1.2/2+1.75}) node[black,below, anchor=south west]
{$\Lcal_2^g, \; g < \frac{\sqrt{q}(\sqrt{q}-1)}{2}$} ;
\draw [densely dashed] ({(1-57^(1/2))/8},0) node[below,yshift=-0.05cm] {$\mu_2$}-- ({(1-57^(1/2))/8},{(1-57^(1/2))/16+3/4}) node[above, xshift=0.2cm,yshift=0.1cm] {$A_2^g$};
\draw [gray,opacity=0.5,->] (-1.4,0) -- (1.4,0) ;
\draw (1.4,0) node[below left] {$x_1$};
\draw [gray,opacity=0.5,->] (0,-1.4) -- (0,2) ;
\draw (0,2) node[below left] {$x_2$};
\draw (-1,0) node[below left] {$-1$}; \draw (1,0) node[below left] {$1$};
\draw (0,-1) node[below left] {$-1$}; \draw (0,1) node[below left] {$1$};
\end{tikzpicture}
\end{center}
\caption{The Weil domain $\overline{\Dcal}_2$. If $X$ is a curve over $\FF_q$ of 
genus $g \geq \frac{\sqrt{q}(\sqrt{q}-1)}{2}$, 
then $P_2(X)=(x_1, x_2)$ as defined by~(\ref{def_x_i}) lies on the grey domain $\overline{\Dcal}_2 \cap \Hcal_2^g$, so that
$x_1$ is lower bounded  by some
constant $\mu_2 = \min \{x_1(P), P \in \Dcal_2\cap \Hcal_2^g\} >-1$, improving Weil bound by Proposition~\ref{T(X)DansDncapHng}.}\label{figure_Ihara_domain}
\end{figure}
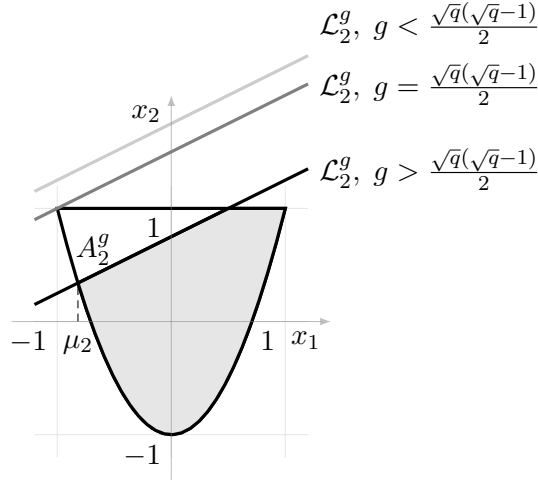

For~$n = 2$,  the domain~$\Dcal_2$ recursively
corresponds by item~\ref{def_Dn_rec} of Proposition~\ref{DnRecursif} and~\eqref{G2}
to the set of~$(x_1,x_2) \in \RR^2$
satisfying:
$$
\begin{cases}
G_1(x_1) = 1-x_1^2 > 0 \\
G_2^+(x_1,x_2) = 1 - x_2 > 0 \\
G_2^-(x_1,x_2) = 1+ x_2 -2x_1^2 > 0
\end{cases},
\qquad \hbox{that is} \qquad
\begin{cases}
-1 < x_1 < 1 \\
2x_1^2 -1 < x_2 < 1
\end{cases}.
$$
We represent the second Weil domain $\overline{\Dcal}_2$ on figure~\ref{figure_Ihara_domain}.

The second Ihara line~$\Lcal_2^g$ with positive slope $\alpha=\frac{1}{\sqrt{q}}$ meets the domain $\overline{\Dcal}_2$ if and only if the genus
is greater than the genus~$g_2$ for which~$\Lcal_2^{g_2}$ contains the point~$(-1,1)$
(see figure~\ref{figure_Ihara_domain}), that is
$$
1 = \alpha\times (-1) + \frac{1}{2g_2\alpha}\left(\frac{1}{\alpha}-\alpha\right)
\qquad\Longleftrightarrow\qquad
g_2 = \frac{\sqrt{q}\left(\sqrt{q}-1\right)}{2}.
$$
For~$g > g_2$, the constraint~$h_2^g(x_1,x_2) \leq 0$ restricts
the domain $\overline{\Dcal}_2$ to the grey domain
$\overline{\Dcal}_2 \cap \Hcal_n^g$, as can be seen on figure~\ref{figure_Ihara_domain}. The point~$A_2^g$ such
that~$[A_2^g,B_2^g] = \overline{\Dcal}_2 \cap \Lcal^g_2$
lies on the curve~$G_2^- = 0$, and minimizes\footnote{It can also be trivially proved
using criteria~\ref{criteria_min_x_1}
since~$\nabla G_2^- = \begin{psmallmatrix}-4x_1\\1\end{psmallmatrix}$
and~$x_1(A_2^g) < 0$.} the first coordinate~$x_1$
on~$\overline{\Dcal}_2 \cap \Hcal_2^g$.

$\mu_2$ is thus the solution in $[-1, 0]$ of the quadratic equation
$G_2^-(x_1, \frac{x_1}{\sqrt{q}}+\frac{q-1}{2g})=0$. Solving it, we find
$$\mu_2= \frac{\alpha^2-\sqrt{\alpha^2+8\left(\frac{q-1}{2g}-1\right)}}{4},$$
so that using~(\ref{def_x_i}),
 we recover the well known Ihara  bound
\begin{equation} \label{Ihara}\tag{\bfseries Ihara  ``second order Weil" bound}
\boxed{
\sharp \Cun(\FF_q) - (q+1)
\leq
\frac{\sqrt{(8q+1)g^2+4q(q-1)g}-g}{2}
}
\end{equation}



\subsection{Weil bounds of higher finite orders} \label{NewBounds}

Following the same line, one can study
Weil bounds {\em of order} $n$ for~$n \geq 3$. We compute in Section~\ref{s_ordre3} the exact
formulae for the Weil bound of order~$3$ and for the genus bound~$g_3$
from which this new bound is better than the Weil bound of order~$2$.

For~$n\geq 4$, computations
to obtain explicit formula for the $n$-order Weil bound becomes intractable.
Therefore we choose to develop an algorithm which, given a genus~$g$
and a size field~$q$, computes the best upper order bound for the number of $\FF_q$-rational points of a curve of genus $g$, together with the corresponding
order $n$. We need to prove  in Section~\ref{GeneralFeatures} some  preliminary results to justify this algorithm.

To illustrate the efficiency of the algorithm we display in Section~\ref{table} a table of numerical
results for orders $n \geq 4$ and few given $q$, $g$.

\subsubsection{Weil bound of order~$3$}\label{s_ordre3}

For~$n=3$, we have by~(\ref{G3}) that
\begin{align*}
G_3^+(x_1,x_2,x_3)
&=
\begin{vmatrix}
1-x_3 & x_1 - x_2\\
x_1 - x_2 & 1-x_1
\end{vmatrix}
&= -(1-x_1)x_3+1-x_1-(x_1-x_2)^2,\\
G_3^-(x_1,x_2,x_3)
&=
\begin{vmatrix}
1+x_3 & x_1 + x_2\\
x_1 + x_2 & 1+x_1
\end{vmatrix}
&= (1+x_1)x_3+(1+x_1)-(x_1+x_2)^2,
\end{align*}
where by~\eqref{G1} we have $G_1^-=1+x_1$ and $G_1^+=1-x_1$. Hence, by item~\ref{def_Dn_graphe} of Proposition~\ref{DnRecursif}, we have
\begin{align*}
\widetilde{G}_3^-
&=-\left(1+x_1-(x_1+x_2)^2\right),\\
\widetilde{G}_3^+
&=1-x_1-(x_1-x_2)^2,
\end{align*}
so that by item~\ref{def_Dn_graphe} of Proposition~\ref{DnRecursif}, we have
$$\Dcal_3 = \left\{(x_1,x_2, x_3) \mid (x_1,x_2) \in \Dcal_2 \quad \hbox{and} \quad
-1 + \frac{(x_1+x_2)^2}{1+x_1}
<
x_3
<
1 - \frac{(x_1-x_2)^2}{1-x_1}\right\}.
$$

The boundary $\partial \Dcal_3 = \overline{\Dcal}_3 \setminus \Dcal_3$ 
is by item~\ref{def_Dn_graphe} of Proposition~\ref{DnRecursif} the part of the 2-dimensional graphs~$G_3^\pm = 0$ above~$\overline{\Dcal}_2$, and both
graphs meet by Lemma~\ref{lem_2G_n} along a curve above the plane 
curve~$\{(x_1, x_2) \in ]-1, +1[^2 ; G_2(x_1, x_2)=0\}$,
where $G_2$ is given by~\eqref{G2}. More precisely,
above the part~$G_2^+ = 0$, that is above the locus~$x_2=1$, the two surfaces~$G_3^- = 0$
and~$G_3^+ = 0$ meet along the segment~$(x_1,1,x_1)$ for $-1\leq x_1 \leq 1$.
Above the curve~$G_2^- = 0$, i.e. above the locus~$x_2 = 2x_1^2-1$,
one has~$x_2+x_1 = (x_1+1)(2x_1-1)$ and~$x_2-x_1 = (x_1-1)(2x_1+1)$.
Therefore, $\{G_3^+=0\}$ and $\{G_3^-=0\}$ above $x_2=2x_1^2-1$ are respectively
\begin{align*}
x_3+1 = \frac{(x_1+x_2)^2}{x_1+1} = (x_1+1)(2x_1-1)^2 = 4x_1^3-3x_1+1 \hbox{ and } x_2=2x_1^2-1\\
x_3-1 = \frac{(x_1-x_2)^2}{x_1-1} = (x_1-1)(2x_1+1)^2 = 4x_1^3-3x_1-1 \hbox{ and } x_2=2x_1^2-1,
\end{align*}
so that their intersection above $x_2=2x_1^2-1$ is the curve~$(x_1,2x_1^2-1,4x_1^3-3x_1)$
for $-1\leq x_1 \leq 1$.

\medskip

\begin{figure}
\begin{center}
\includegraphics[scale=0.5]{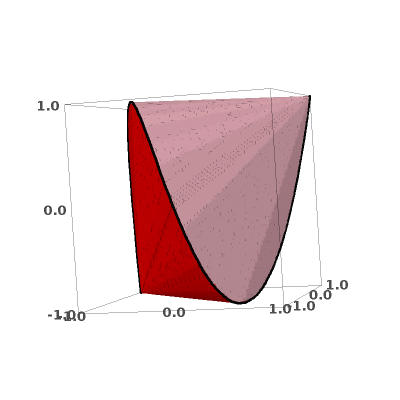}
\includegraphics[scale=0.5]{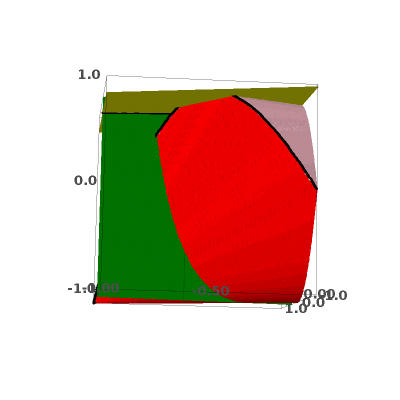}
\includegraphics[scale=0.5]{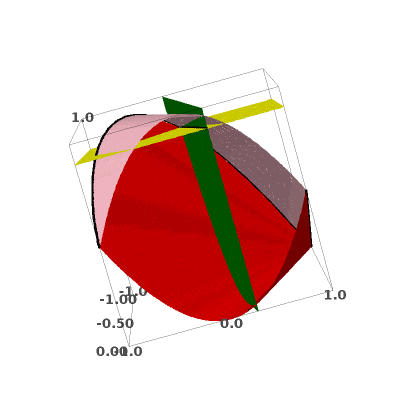}
\includegraphics[scale=0.5]{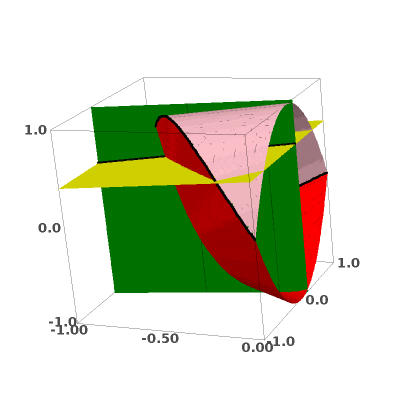}
\end{center}
\caption{a complete view of the third Weil domain~$\overline{\Dcal}_3$ (upper left) and three partial views of $\overline{\Dcal}_3\cap \Lcal_3^g$ (for $-1 \leq x_1 \leq 0$): one for $g<g_3$ (upper right), one for $g=g_3$ (bottom left) and one for $g>g_3$ (bottom right). The red (resp. pink) surface is the locus $\{G_3^-=0\}$ (resp. $\{G_3^+=0\}$) above $\overline{\Dcal}_2$. The dark green almost vertical plane is $x_2=\alpha x_1+\frac{q-1}{2g}$, the light green almost horizontal plane is $x_3=\alpha^2 x_1+\frac{q^2-1}{2g\sqrt{q}}$. Their intersection is the black third Ihara line $\Lcal_3^g$ in genus $g$, which cuts the red surface only for $g \geq g_3$.}
\label{figure_Weil_3}
\end{figure}

\begin{itemize}
\item Looking at figure~\ref{figure_Weil_3}, for~$g$ small the line~$\Lcal_3^g$ does not meet the third Weil domain~$\overline{\Dcal}_3$.
\item As $g$ increase, it intersects the domain inside the~$G_3^+ = 0$ part,
and increasing again inside the~$G_3^- = 0$ part.
This happens for~$g$ greater than the value~$g_3$ such that the
line~$\Lcal_3^{g_3}$ crosses the boundary~$\partial \Dcal_3$ where the two
graphs~$G_3^\pm = 0$ above $\overline{\Dcal}_2$ meet. In other terms, the value~$g_3$
is such that there exist a point
\begin{equation} \label{P3g3DansTout}P_3^g(x_1) \in \Lcal_3^{g_3} \cap \overline{\Dcal}_3 \cap \{G_3^-=0\} \cap \{G_3^+=0\}.
\end{equation}

We now compute this value $g_3$, and the corresponding point~$P_3^{g_3}(x_1) \in \Lcal_3^{g_3} \cap \overline{\Dcal}_3 \cap \{G_3^-=0\} \cap \{G_3^+=0\}$.
Using~\eqref{EquationsToits} and~\eqref{defLng}, there exists some $x_1 \in \RR$, such that $(x_1, \frac{1}{g_3})$ is a solution of the two polynomials in $(x_1, \frac{1}{g})$
\begin{equation}\label{ResG3}
G_3^-
\left(x_1, \frac{x_1}{\sqrt{q}}+\frac{q-1}{2g},\frac{x_1}{q}+\frac{q^2-1}{2g\sqrt{q}}\right)
=0
\quad\text{and}\quad
G_3^+
\left(x_1, \frac{x_1}{\sqrt{q}}+\frac{q-1}{2g},\frac{x_1}{q}+\frac{q^2-1}{2g\sqrt{q}}\right)
=0.
\end{equation}
Eliminating $x_1$ in these equations, $\frac{1}{g_3}$ is a
solution of their resultant in~$x_1$. We have avoided calculations and factorization by hand: with the help of {\tt magma},
we obtain that~$g_3$ is a root of
$$
(g-1)
(g-q)
\left(g-\frac{\sqrt{q}(q-1)}{\sqrt{2}}\right)
\left(g+\frac{\sqrt{q}(q-1)}{\sqrt{2}}\right).
$$
Of course, $g= -\frac{\sqrt{q}(q-1)}{\sqrt{2}} <0$ cannot be a solution coming from a curve $X$.
If $g=1$ or~$q$, substituting these values for $g$ in~\eqref{ResG3}, we find that~$x_1= \frac{q+1}{2\sqrt{q}}$
is the only corresponding solution (for both values $g=1$ and $g = q$),
which is not contained in~$\left[-1, +1\right]$, so that these two values for the genus cannot come from a curve over a finite field
by Weil bound (of order $1$).
Therefore
$$
g_3=\frac{\sqrt{q}(q-1)}{\sqrt{2}},
$$
and substituting again this value in~\eqref{ResG3}, we find (using {\tt magma}) that $x_1=-\frac{\sqrt{2}}{2}$ is the only corresponding solution, so that, as seen above,  $x_2=2x_1^2-1=0$ and $x_3=4x_1^3-3x_1=\frac{\sqrt{2}}{2}$. Finally, $g=g_3=\frac{\sqrt{q}(q-1)}{\sqrt{2}}$ is the only solution, corresponding to the point $P_3(-\frac{\sqrt{2}}{2})=(-\frac{\sqrt{2}}{2}, 0, \frac{\sqrt{2}}{2})$, for which~\eqref{P3g3DansTout} holds.

\item For~$g > g_3$, the
intersection of the line~$\Lcal_3^g$ with the domain~$\overline{\Dcal}_3$
is a segment~$[A_3^g,B_3^g]$,
where the point~$A_3^g$ lies on~$G_3^- = 0$. 
\end{itemize}


The point of the proof of the next Theorem is to check that, for $g > g_3$, the point $P_0=A_3^g$ satisfies the requirements of criteria~\ref{criteria_min_x_1}, showing that
this point minimizes~$x_1$ on the domain~$\overline{\Dcal}_3\cap \Hcal_3^g$, which contains all points coming from curves by Proposition~\ref{T(X)DansDncapHng}.

\begin{theo} \label{Ordre3}
Let~$X$ be a smooth projective absolutely irreducible curve over~$\FF_q$
of genus~$g$ such that~$g > g_3=\frac{\sqrt{q}(q-1)}{\sqrt{2}}$. Then, we have
$$\sharp \Cun(\FF_q)-q-1
\leq
\Biggl(
\sqrt{a(q)+\frac{b(q)}{g}+\frac{c(q)}{g^2}} - 1+ \frac{1}{q}+\frac{d(q)}{g} \Biggl)g\sqrt{q},$$
where
\begin{equation*}
\left\{
\begin{matrix}
a(q)&=& 5+\frac{8}{\sqrt{q}} -\frac{1}{q^2}\\
b(q)&=& \frac{q-1}{q\sqrt{q}}\Bigl(q^2-4q \sqrt{q}+2q+4\sqrt{q}-1\Bigr)\\
c(q)&=& \frac{q-1}{4 q}\Bigl(q^3-5q^2-8q\sqrt{q}-5q-8\sqrt{q}+1\Bigr)\\
d(q) &=& \frac{2\sqrt{q}(q-1)^2}{q}.
\end{matrix}
\right.
\end{equation*}
\end{theo}

\begin{preuve}
Let $g>g_3$. Referring to the last item in the discussion above the statement of the Theorem,  we begin to prove that the unique $-1<\mu_3<0$, such that $A_3^g=P_3(\mu_3) \in \overline{\Dcal}_3\cap \Lcal_3^g\cap \{G_3^-=0\}$, satisfies criteria~\ref{criteria_min_x_1}. Then, we compute the exact value of $\mu_3$, which is, by criteria~\ref{criteria_min_x_1} together with Proposition~\ref{T(X)DansDncapHng}, a lower bound for $x_1$'s coming from curve over $\FF_q$. We deduce the Theorem using~\eqref{def_x_i}.

\bigskip

Let us prove that $P_0=P_3(\mu_3))$ satisfies the requirements of criteria~\ref{criteria_min_x_1}.
The first requirement holds true: $A_3^g\in \{G_3^-=0\}$. We prove later in corollary~\ref{DeriveeDirectionnelle} that the third one also holds true for any $n\geq 1$. we have only to check the second requirement, that
$\partial_i G_3^-(P_3^g(\mu_3)) \geq 0$ for $i=2, 3$. Since by~\eqref{G3} we have
$$G_3^-(x_1,x_2,x_3)
=
\begin{vmatrix}
1+x_3 & x_1 + x_2\\
x_1 + x_2 & 1+x_1
\end{vmatrix},$$
we deduce that for any $(x_1, x_2, x_3) \in \RR^3$,
$$\left\{
\begin{matrix}
\partial_2 G_3^-(P_3(x_1)) &= -2(x_1+x_2)  \\
\partial_3 G_3^-(P_3(x_1)) &= 1+x_3 .
\end{matrix}
\right.
$$
If $x_1=\mu_3$, then $P_0=A_3^g=P_3(\mu_3) \in \overline{\Dcal}_3 \subset ]-1, +1[^3$, hence
$$\partial_3G_3^-(P_3(\mu_3)) = 1+x_3(P_0)>0,$$
so that it remains to prove that $x_1(P_0)+x_2(P_0)\leq 0$. 

Since $P_3(\mu_3)=A_3^g\in\{G_3^-=0\}\cap\Lcal_3^g$, we have
$$
\left\{
\begin{matrix}
(x_1+x_2)^2 &= &(1+x_1)(1+x_3)\\
x_3 &= &\alpha^2x_1+\frac{q^2-1}{2g\sqrt{q}}.
\end{matrix}
\right.
$$
We deduce that $P_3(\mu_3)$ lies on on the plane parabola
$$x_1+x_2 = \pm \sqrt{(1+x_1)\left(1+ \alpha^2x_1+\frac{q^2-1}{2g\sqrt{q}}\right)},$$
and are reduced to check that  $P_3(\mu_3)$ lies on on the lower branch
$$x_1+x_2 = - \sqrt{(1+x_1)\left(1+ \alpha^2x_1+\frac{q^2-1}{2g\sqrt{q}}\right)}.$$ 
The point is that $(x_1, x_2)=(-1, 1)$ lies on on both branches of this parabola, and that the tangent line to the parabola at this point is the vertical one. Hence, the upper branch defines a concave subset
$$\left\{(x_1, x_2) \mid x_1 \geq -1, \quad  x_2\geq1 \quad \hbox{and} \quad x_1+x_2 \leq \sqrt{(1+x_1)\left(1+ \alpha^2x_1+\frac{q^2-1}{2g\sqrt{q}}\right)}\right\}$$
of $\RR^2$, containing $(-1, 1)$ and $(0, \sqrt{1+\frac{q^2-1}{2g\sqrt{q}}})$ with $\sqrt{1+\frac{q^2-1}{2g\sqrt{q}}}>1$. It follows that for any $-1 < x_1 < 0$, a point $(x_1, x_2)$ lying on the upper parabolic branch lies also on the locus $x_2>1$, which cannot hold
for $P_0=P_3(\mu_3)) \in \overline{\Dcal}_3 \subset [-1, -1]^3$. 

\medskip

It follows that all requirements of criteria~\ref{criteria_min_x_1} hold true for the point $P_3(\mu_3)$, which lies on the branch $x_1+x_2 = - \sqrt{(1+x_1)(1+ \alpha^2x_1+\frac{q^2-1}{2g\sqrt{q}})}$, so that $\mu_3$ is a lower bound on $\overline{\Dcal}_3\cap \Hcal_3^g$ for $g >g_3$.

Since $P_3(\mu_3)$ lies on $\Lcal_3^g$, it lies in particular on $x_2=\alpha x_1+\frac{q-1}{2g}$, so that $\mu_3$ is a solution of the quadratic equation
$$(x_1+ \alpha x_1+\frac{q-1}{2g})^2=(1+x_1)(1+ \alpha^2x_1+\frac{q^2-1}{2g\sqrt{q}}),$$
where we recall that $\alpha = \frac{1}{\sqrt{q}}$.
Solving this equation proves that the only solution in $[-1, 0]$ is
\begin{equation*}
\mu_3= \frac{
 \frac{q-1}{q}-\frac{2\sqrt{q}(q-1)^2}{gq} -\sqrt{a(q)+\frac{b(q)}{g}+\frac{c(q)}{g^2}}}{2},
\end{equation*}
hence the Theorem using~(\ref{def_x_i}).
\end{preuve}

The formula becomes nicer if we let~$g$ going to infinity. We obtain the following third
order asymptotic bound for Ihara's constant $A(q)$ (see~\cite{Ihara}).

\begin{cor} The Ihara constant~$A(q)$ is bounded above by:
$$
A(q) \leq \Biggl(\sqrt{5+\frac{8}{\sqrt{q}}-\frac{1}{q^2}} - 1 +\frac{1}{q}\Biggl)\sqrt{q}.
$$
\end{cor}

\begin{rk}
For $q$ large the preceding upper bound is equivalent
to~$(\sqrt{5}-1)\sqrt{q} \simeq 1,236\sqrt{q}$. This is 
better than the upper bound of~$A(q)$ following from Ihara
bound
$$
A(q) \leq \frac{\sqrt{8q+1}-1}{2},
$$
which is equivalent to $\sqrt{2q} \simeq 1,414 \sqrt{q}$
for~$q$ large. We prove later in Theorem~\ref{TVameliore} that we can also recover Dinfeld-Vl{\u{a}}du{\c{t}} and Tsfasman bounds with this viewpoint.
\end{rk}


\subsubsection{General features} \label{GeneralFeatures}

The idea of the algorithm producing numerical values for higher order Weil bounds for given $g, q$ is
the following. It computes
numerically in a first step an approximation of the unique $-1\leq x_1<0$ for which $P_n(x_1)$ lies 
on the intersection of~$\Lcal_n^g$
with~$G_n^- = 0$, and it checks numerically in a second step
that this point satisfies criteria~\ref{criteria_min_x_1}. This algorithm is valid
thanks to Proposition~\ref{Suite_g_n} that such an $x_1$ as in the first step
do exists, at least for~$g$ large enough.

\begin{prop}\label{Suite_g_n}
There exists a sequence~$(g_n)_{n\geq 2}$ such that
~$g > g_n$
if, and only if, the intersection of the $n$-th Ihara line~$\Lcal_n^g$ in genus $g$ with the interior of the $n$-the Weil domain~$\Dcal_n$
is a non-empty $x_1$-segment~$\left]A_n^g,B_n^g\right[$, 
with~$A_n^g$ lies on on the hypersurface of~$\overline{\Dcal}_{n-1} \times \RR$
having equation~$G_n^- = 0$, and with $-1\leq x_1<0$.
\end{prop}

\begin{preuve}
We keep notations of Section~\ref{SubsecCurveLocus} and we begin by studying
the infinite line~$
\Lcal_n^\infty
=
\{
(x_1,\alpha x_1, \ldots, \alpha^{n-1}x_1),
\, x_1 \in \RR\}$, whose points are the~$P_n^\infty(x_1)$ for~$x_1\in\RR$.
This line~$\Lcal_n^\infty$ and the domain~$\Dcal_n$ have a non empty
intersection since the origin of~$\RR^n$  belongs to~$\Dcal_n\cap \Lcal_n^\infty$. Since $\Dcal_n$ is convex and open,~$\Dcal_n \cap \Lcal_n^\infty$
must be a non empty segment~$\left]P_n^\infty(\mu_n^{\infty}), P_n^\infty(\nu_n^{\infty})\right[$
for some~$-1\leq \mu_n^{\infty}<\nu_n^{\infty}\leq1$.
Necessarily~$\mu_n^{\infty} < 0$
and\footnote{If~$\Ocal$ is a convex open subspace of~$\RR^n$, for any~$x\in\Ocal$ and
any~$y\in\overline{\Ocal}$, the segment~$\left]x,y\right[$ is included
in~$\Ocal$.}~$\overline{\Dcal}_n \cap \Lcal_n^\infty
=
\left[P_n^\infty(\mu_n^{\infty}),P_n^\infty(\nu_n^{\infty})\right]$.

Next, we prove that the point~$P_n^\infty(\mu_n^{\infty})$ lies on the hypersurface
of~$\Dcal_{n-1}\times \RR$ having equation~$G_n^-=0$.
To relieve the notations, for~$\epsilon=\pm$ or nothing, we
put~$L_n^\epsilon(x_1) = G_n^\epsilon(x_1,\alpha x_1, \ldots, \alpha^{n-1}x_1)$.
Let us prove by induction that the abscissa~$\mu_n^{\infty}$ satisfy
\begin{equation}\label{eq_hyp_rec}
\forall n\geq 1,
\qquad
L_n^-(\mu_n^{\infty}) = 0
\text{ and }
L_m(\mu_n^{\infty}) > 0, \; \forall m < n.
\end{equation}
For~$n=1$ this is true from~\eqref{G1}:~$\mu_1 = -1$ is the only root of~$L_1(x_1) = 1+x_1$
and satisfies~$L_0(-1) > 0$ since~$L_0 = 1$. Suppose that assertion~(\ref{eq_hyp_rec}) holds true for some~$n\geq 1$. By Lemma~\ref{lem_2G_n}, one has
\begin{equation*}\label{eq_2L_n}
2\underbrace{L_n(\mu_n^{\infty})}_{=0}
=
\underbrace{L_{n-1}^+(\mu_n^{\infty})}_{>0}L_{n+1}^-(\mu_n^{\infty})
+
\underbrace{L_{n-1}^-(\mu_n^{\infty})}_{>0}L_{n+1}^+(\mu_n^{\infty}),
\end{equation*}
where the underbrace (in)equalities come from the induction hypotheses.
We first establish that~$L_{n+1}^+(\mu_n^{\infty}) > 0$. To this end,
let~$m$ be the quotient of~$n+2$ by~$2$. By Lemma~\ref{lem_factorization},
one has~$G_{n+1}^+ = \Det(T_{m} - H_m)$, where~$T_m$ is a Toeplitz given by~\eqref{Tn+1}, while~$H_m$
is Hankel given by~\eqref{Hn+1}. Specializing
to the point~$(\mu_n^{\infty},\alpha \mu_n^{\infty}, \ldots, \alpha^{n} \mu_n^{\infty})$, the Hankel part
becomes a rank one matrix: 
$$
-H_m
=
\alpha^{\frac{(-1)^{n+1} + 1}{2}}
\times
(-\mu_n^{\infty})
\times
{}^t(\alpha^{m-1},\ldots,\alpha,1)
\times
(\alpha^{m-1},\ldots,\alpha,1).
$$
This matrix is clearly symmetric positive semi definite.
But, for general~$A$ and~$B$ positive semi-definite symmetric
matrices, one has~$\det(A+B) \geq \det(A)$
(easy consequence of \cite[Corollary~4.3.3]{MatrixAnalysis}).
We deduce by induction hypothesis that, for some $z\in \RR^{2n}$,
\begin{align*}
L_{n+1}^+(\mu_n^{\infty})
&=
\Det\left(
T_m(1,\mu_n^{\infty},\ldots,\alpha^{m-2}\mu_n^{\infty})
+
H_m(z)
\right) \\
&\geq
\Det(T_m(1,\mu_n^{\infty},\ldots,\alpha^{m-2}\mu_n^{\infty}))
=
G_{m-1}(\mu_n^{\infty},\ldots,\alpha^{m-2}\mu_n^{\infty})
=
L_{m-1}(\mu_n^{\infty}) > 0
\end{align*}
holds true. It then follows by~(\ref{eq_2L_n}) that~$L_{n+1}^-(\mu_n^{\infty}) <0$.
Since~$L_{n+1}(0) = 1$, there exists~$\mu \in \left]\mu_n^{\infty},0\right[$ such
that~$L_{n+1}(\mu) = 0$. Since the points~$P_n^\infty(x_1)$
with~$x_1 \in \left]\mu_n^{\infty},0\right[$ of the line~$\Lcal_n^\infty$ belong
to~$\Dcal_n$, this is the case for the point~$P_n^\infty(\mu)$.
Hence one must have~$L_m(\mu) > 0$ for all~$m \leq n$ and~$L_{n+1}(\mu) = 0$.
This proves that point~$P_{n+1}^\infty(\mu)$ belongs
to~$\overline{\Dcal}_{n+1}$ and that~$\mu = \mu_{n+1}$. This concludes
the induction, hence in particular $G_n^-(P_n^{\infty}(\mu_n^{\infty}))=L_n^-(P_n^{\infty}(\mu_n^{\infty}))=0$
as asserted, with moreover $-1\leq \mu_n^{\infty} < 0$.

Now for finite genus, formula~(\ref{EquationsToits}) shows that the union ${\mathcal P}_n = \cup_{g \in \left]0, +\infty \right]} \Lcal_n$ is a closed half linear $2$-dimensional plane
with parameters $\left(x_1 \in \RR, \frac{1}{g} \in \left[0, +\infty\right[\right)$. Hence the intersection ${\mathcal P}_n \cap \Dcal_n$ is
a $2$-dimensional convex subset in the closed half plane ${\mathcal P}_n$. Moreover, we have just seen above that the cut out with the line having parameter $\frac{1}{g}=0$ meet it on an non-empty segment with an end $A_n^{\infty}=P_n(\mu_n^{\infty}) \in \{G_n^-=0\}$ with $-1 \|eq \mu_n^{\infty} < 0$. Since by Proposition~\ref{ConvexeConcave} the locus $\{(x_1, \ldots, x_n) \in \Dcal_{n-1} \mid G_n^-(x_1, \ldots, x_n)>0\}$ is convex, this segment is non empty and meet $G_n^-=0$ for at least one point with $x_1<0$ for a connected subset of the parameter $0 \leq \frac{1}{g} < +\infty$ containing $0$, hence have the form $0\leq \frac{1}{g} < \frac{1}{g_n}$ for some $0 < g_n <+\infty$. for such a $g$, the cut out with $\Lcal_n^g$ is a segment $[A_n^g, B_n^g]$, with $A_n^g\neq B_n^g$ since $g>g_n$, with $A_n^g \in \{G_3^-=0\}$ and with $-1 \leq x_1(A_n^g)<0$. The Proposition is proved.
\end{preuve}

It turns out that under the assumption that a point on $\Lcal_n^g$ lies on the locus
$\overline{\Dcal}_n \cap \{G_n^-=0\}$ always satisfy the last requirement of criteria~\ref{criteria_min_x_1}:

\begin{cor} \label{DeriveeDirectionnelle}
Let $n \in \NN$, $g>g_n$ and $\mu_n \in ]-1, 0[$ such that $P_n^g(\mu_n)$ is the point $A_n^g$ of Proposition~\ref{Suite_g_n}. Then
$\sum_{i=1}^n \partial_i G_n^-(P_n^g(\mu_n)) \alpha^{i-1} > 0$.
\end{cor}
\begin{preuve}
Denote by $\varphi(x_1) = G_n^-(P_n^g(x_1))$ for $x_1 \in {\mathbb R}$. 
This is a ${\mathcal C}^{\infty}$ function of the real variable $x_1$.

By assumption, $P_n^g(\mu_n) \in \{G_n^-=0\}\cap \Lcal_n^g$ and $g > g_n$, hence by Proposition~\ref{DnRecursif} and Proposition~\ref{Suite_g_n}
we have $P_n^g(\mu_n)=A_n^g$ and $]P_n^g(\mu_n), B_n^g[ \subset \Dcal_n \subset \{G_n^->0\}$, and there exists some $\varepsilon>0$, such that
\begin{equation} \label{CaCraindra}
\forall x_1 \in \RR, \quad \mu_n < x_1 < \mu_n + \varepsilon \Longrightarrow P_n^g(x_1) \in \{G_n^->0\}.
\end{equation}
Hence by item~\ref{def_Dn_rec} of Proposition~\ref{DnRecursif}, we have $P_n^g(x_1) \in \{G_n^->0\}$ for $\mu_n < x_1 < \mu_n+\varepsilon$, so that $\varphi(x_1)>0$ for such $x_1$. Since
by assumption $\varphi(\mu_n)=0$, we deduce that
\begin{equation*} \label{DeriveePhi}
0 \leq \frac{\mathrm{d} \varphi}{\mathrm{d} t}(\mu_n) = \frac{\mathrm{d}  G_n^-(P_n^g(\mu_n))}{\mathrm{d} t}
= \sum_{i=1}^n \partial_i G_n^-(P_n^g(\mu_n)) \alpha^{i-1}=\scalaire{\nabla G_n^-(\mu_n)}{\vec{u}},
\end{equation*}
where $\vec{u}=(1, \alpha, \ldots, \alpha^{n-1})$ is a director vector of $\Lcal_n^g$.
Now, suppose by contradiction that $ \frac{\mathrm{d} \varphi}{\mathrm{d} t}(\mu_n) =0$. Then by~\eqref{DeriveePhi} the line $\Lcal_n^g$ would be perpendicular to $\nabla G_n^-(\mu_n^{\infty})$, hence
tangent to $\{G_n^-=0\}$ at $P_n^g(\mu_n)=A_n^g$ by Proposition~\ref{Suite_g_n}. But
$\{G_n^-\geq 0\}\cap \overline{\Dcal}_{n-1}$ is convex by Proposition~\ref{ConvexeConcave}, so that $\Lcal_n^g$
would never meet the interior $\{G_n^-=0\}$, a contradiction with~\eqref{CaCraindra}.
\end{preuve}

It turns out that proving that this points is also always a minimum of the
coordinate function~$x_1$ seems to be difficult for $N \geq 4$. But this
is an easy task numerically using criteria~\ref{criteria_min_x_1}:
this will be done in our {\tt magma} routine in Section \ref{table}.

\subsubsection{Numerical results for $n \geq 4$} \label{table}

We have
implemented a routine in {\tt magma} and leading to the following
results.

In the following tabular we give, for genus~$g$ such that~$1 \leq g\leq 52$,
and for~$q \in \{2,3\}$, the best generalized Weil bound for the number
of points on a curve of genus~$g$ over~$\FF_q$.
The number in brackets corresponds to the {\em optimal} order Weil bound.
For instance, if
this number is~$2$, this means that the bound is the well known Ihara bound
and that bounds of other orders are not better.
$$
\begin{array}{|*{4}{ccc|}}
\hline
g& q=2& q=3&g& q=2& q=3&g& q=2&q=3&g& q=2&q=3\\*
\hline
1& 5(3)& 7(2)& 14& 16(6)& 26(4)& 27& 25(8)& 42(5)& 40& 34(9)& 57(6) \\
2& 6(3)& 9(2)& 15& 17(7)& 28(4)& 28& 26(8)& 43(5)& 41& 35(9)& 58(6) \\
3& 7(4)& 10(3)& 16& 18(7)& 29(5)& 29& 27(8)& 44(5)& 42& 35(9)& 59(6) \\
4& 8(4)& 12(3)& 17& 18(7)& 30(5)& 30& 27(8)& 46(5)& 43& 36(9)& 60(6) \\
5& 9(5)& 14(3)& 18& 19(7)& 31(5)& 31& 28(8)& 47(5)& 44& 37(9)& 61(6) \\
6& 10(5)& 15(3)& 19& 20(7)& 32(5)& 32& 29(8)& 48(6)& 45& 37(9)& 62(6) \\
7& 11(5)& 17(4)& 20& 21(7)& 34(5)& 33& 29(8)& 49(6)& 46& 38(9)& 63(6) \\
8& 11(5)& 18(4)& 21& 21(7)& 35(5)& 34& 30(8)& 50(6)& 47& 38(9)& 65(6) \\
9& 12(6)& 19(4)& 22& 22(7)& 36(5)& 35& 31(8)& 51(6)& 48& 39(9)& 66(6) \\
10& 13(6)& 21(4)& 23& 23(7)& 37(5)& 36& 31(8)& 52(6)& 49& 40(9)& 67(6) \\
11& 14(6)& 22(4)& 24& 23(8)& 38(5)& 37& 32(8)& 54(6)& 50& 40(9)& 68(6) \\
12& 15(6)& 24(4)& 25& 24(8)& 40(5)& 38& 33(9)& 55(6)& 51& 41(9)& 69(6) \\
13& 15(6)& 25(4)& 26& 25(8)& 41(5)& 39& 33(9)& 56(6)& 52& 42(9)& 70(6) \\
\hline
\end{array}
$$
Comparing with the results available on the web site
{\tt http://www.manypoints.org/}, we unfortunately observe that we do not
beat any record, and that we reach the records exactly in those cases were it is held by Osterl\'e bounds! . So we are pretty sure that
there is a link
between Osterle bound and the one in this Section, even if we do not know how to relate
the two approaches.

\subsection{Drinfeld-Vl{\u{a}}du{\c{t}} and Tsfasman bounds as a bound of infinite order} \label{asymptotique}

In this Section, we recover the 
Drinfeld-Vl{\u{a}}du{\c{t}} and Tsfasman bounds for the
asymptotic of the number of points on curves over~$\FF_q$, giving a new meaning for the defect $\delta$ defined in~\cite{T92}. With our point of
view, these bounds can be considered as Weil bounds of infinite order.

We consider a sequence~$(X_n)_{n \geq 1}$ of absolutely irreducible smooth projective curves
over~$\FF_q$. Let~$(g_n)_{n\geq 1}$ be the genus sequence
and for~$r\geq 1$ and let~$B_r(X_n)$ be the number of points of degree~$r$
on~$X_n$:
$$
B_r(X_n)
=
\sharp
\left\{
P \in X_n(\overline{\FF}_q) \mid \deg(P) = r
\right\}.
$$
Following Tsfasman, we assume this sequence
to be {\em asymptotically exact}: this means
that~$\lim_{n\to+\infty} g_n = +\infty$
and that for any~$r\geq 1$, the
sequence~$\left(\frac{B_r(X_n)}{g_n}\right)_{n\geq 1}$
admits a limit. Then, as usual, we put
$$
\forall r\geq 1,
\qquad
\beta_r
\overset{\text{def}}{=}
\lim_{n\to+\infty}
\frac{B_r(X_n)}{g_n}.
$$

In the following Theorem, the notations and the normalizations are the same
as in Definitions~\ref{def_scalar_product} and~\ref{Normalization}.

\begin{theo} \label{TVameliore}
Let $(\Cun_n)_{n \geq 1}$ be a sequence of absolutely irreducible smooth projective curves
over~$\FF_q$.
If this sequence is asymptotically exact, then its {\em defect} $\delta$, defined as
$\delta = 1-\sum_{r=1}^{\infty} \frac{r\beta_r}{\sqrt{q^r}-1}$, satisfy
$$
\delta
=
\lim_{m \to +\infty} \lim_{n \to +\infty}
\frac{
\left\Vert
\gamma_{X_n}^0+ \gamma_{X_n}^1+\cdots+\gamma_{X_n}^{m-1}
\right\Vert_{X_n}^2
}{m}.
$$
\end{theo}

\begin{preuve}
Thanks to the
Gram matrix following Definition~\ref{Normalization}, we compute
\begin{align*}
\frac{1}{m}
\left\Vert \sum_{i=0}^{m-1} \gamma_{X_n}^i \right\Vert_{X_n}^2
&=
\frac{1}{m}
\begin{pmatrix}1&\cdots&1\end{pmatrix}
\times
\Gram\left(\gamma_{X_n}^0,\ldots,\gamma_{X_n}^{m-1}\right)
\times
\begin{pmatrix}1\\\vdots\\1\end{pmatrix} \\
&=
\frac{1}{m}
\left[
m + 2\sum_{i=1}^{m-1} (m-i) x_i
\right] \\
&=
1 + \frac{1}{m}\sum_{i=1}^{m-1} (m-i)\frac{(q^i+1) - \sharp X_n(\FF_{q^i})}{g_n q^{\frac{i}{2}}} \\
&=
1 + \frac{1}{g_n} \sum_{i=1}^{m-1} (m-i)\left(q^{\frac{i}{2}} + q^{-\frac{i}{2}}\right)
-
\sum_{i=1}^{m-1} \left(1 - \frac{i}{m}\right)
\sum_{r \mid i} \frac{r}{q^{\frac{i}{2}}}\times\frac{B_r(X_n)}{g_n}.
\end{align*}
Letting~$n$ tends to~$+\infty$ leads to
\begin{equation}\label{eq_n_to_infty}
\lim_{n\to+\infty}
\frac{1}{m}
\left\Vert \sum_{i=0}^{m-1} \gamma_{X_n}^i \right\Vert_{X_n}^2
=
1 - \sum_{rs\leq m-1} \left(1 - \frac{rs}{m}\right)
\frac{r\beta_r}{q^{\frac{rs}{2}}}
=
1 - \sum_{r=1}^{m-1}
\left(
\sum_{s=1}^{\lfloor\frac{m-1}{r}\rfloor} 
\left(1 - \frac{rs}{m}\right)
\frac{1}{q^{\frac{rs}{2}}}
\right)
r \beta_r.
\end{equation}
To conclude, for any~$r,m\geq 1$,we remark that
$$
\frac{1}{q^{\frac{r}{2}} -1}
-
\sum_{s=1}^{\lfloor\frac{m-1}{r}\rfloor} 
\left(1 - \frac{rs}{m}\right)
\frac{1}{q^{\frac{rs}{2}}}
=
\sum_{s=1}^{+\infty} \frac{1}{q^{\frac{rs}{2}}}
-
\sum_{s=1}^{\lfloor\frac{m-1}{r}\rfloor} 
\left(1 - \frac{rs}{m}\right)
\frac{1}{q^{\frac{rs}{2}}}
=
\sum_{s=\frac{m-1}{r}}^{+\infty} \frac{1}{q^{\frac{rs}{2}}}
+
\frac{1}{m}
\sum_{s=1}^{\lfloor\frac{m-1}{r}\rfloor}
\frac{rs}{q^{\frac{rs}{2}}}.
$$
Being the remainder of a converging serie, the first term of the last
expression tends to zero when~$m\to+\infty$; as for the second term it
also tends to zero by Cesaro. Therefore
$$
\lim_{m\to+\infty}
\sum_{s=1}^{\lfloor\frac{m-1}{r}\rfloor} 
\left(1 - \frac{rs}{m}\right)
\frac{1}{q^{\frac{rs}{2}}}
=
\frac{1}{q^{\frac{r}{2}} -1},
$$
and the Theorem follows letting~$m$ tends to~$+\infty$ in~(\ref{eq_n_to_infty}).
\end{preuve}

Just pointing out that a norm is non-negative, Theorem~\ref{TVameliore} leads to the well-known bound:
\begin{equation} \label{TV}\tag{\bfseries Tsfasman}
\boxed{
\sum_{r=1}^{\infty} \frac{r \beta_r}{\sqrt{q^r}-1} \leq 1.
}
\end{equation}
As is well known, Tsfasman bound is itself a refinement of
\begin{equation} \label{DV}\tag{\bfseries Drinfeld-Vl{\u{a}}du{\c{t}}}
\boxed{
\limsup_{n \to +\infty} \frac{\sharp \Cun_n(\FF_q)}{g_n} \leq \sqrt{q} -1.
}
\end{equation}

\section{Relative bounds}\label{s_Relative_bounds}

From now on, we concentrate on the relative situation. Our starting point
is a finite morphism~$f : X \to Y$ of degree~$d$, where~$X$ and~$Y$ are
 absolutely irreducible smooth projective curves defined over~$\FF_q$,
whose genus are denoted by~$g_X$ and~$g_Y$.

\subsection{The relative Neron-Severi subspace}
The morphism~$f\times f$ from~$X\times X$ to~$Y\times Y$ induces two
morphisms
$$
(f\times f)_* : \NS(X\times X)_\RR \longrightarrow \NS(Y\times Y)_\RR
\quad\text{and}\quad
(f\times f)^* : \NS(Y\times Y)_\RR \longrightarrow \NS(X\times X)_\RR
$$
satisfying~$(f\times f)_* \circ (f\times f)^* = d^2\Id_{Y\times Y}$. We put:
\begin{equation} \label{phipsi}
\varphi = \frac{1}{d} (f\times f)_*
\qquad\text{and}\qquad
\psi = \frac{1}{d} (f\times f)^*.
\end{equation}
Each of these morphisms can be restricted to the  euclidean 
spaces~$(\HVperp_X, \scalaire{\cdot}{\cdot}_X)$
and~$(\HVperp_Y, \scalaire{\cdot}{\cdot}_Y)$
associated to the curves~$X$ and~$Y$ as in Section~\ref{s_Euclidean_space}.
Recall that both are the image
of the orthogonal projections given by~\eqref{projection}:
$$
p_X : \NS(X\times X)_\RR \longrightarrow \HVperp_X
\qquad\text{and}\qquad
p_Y : \NS(Y\times Y)_\RR \longrightarrow \HVperp_Y.
$$
In fact, we still restrict a little bit more the morphisms~$\varphi$
and~$\psi$. Let~$F_X : X \to X$ and~$F_Y : Y \to Y$ be the $q$-Frobenius
on~$X$ and~$Y$. As usual, for~$k\geq 0$, we denote by~$\Gamma_{F_X}^k$
 the class in~$\NS(X\times X)_\RR$ of the $k$-th iterated of~$F_X$.
We do the same for~$\Gamma_{F_Y}^k$ inside~$\NS(Y\times Y)_\RR$.

\begin{defi} \label{def_gamma_s}
For~$i \geq 0$, we put:
$$
\gamma_X^i = \frac{1}{\sqrt{2q^i}} p_X\left(\Gamma_{F_X}^i\right)
\in \HVperp_X,
\qquad\qquad
\gamma_Y^i = \frac{1}{\sqrt{2q^i}} p_Y\left(\Gamma_{F_Y}^i\right)
\in \HVperp_Y,
$$
and we denote by~$\Frob_X$ and~$\Frob_Y$ the two subspaces of~$\HVperp_X$
and~$\HVperp_Y$ defined by
$$
\Frob_X
=
\Vect
\left(
\gamma_X^i, \; i\geq 0
\right)
\subset
\HVperp_X,
\qquad\qquad
\Frob_Y
=
\Vect
\left(
\gamma_Y^i, \; i\geq 0
\right)
\subset
\HVperp_Y.
$$
\end{defi}

\begin{rk}
Note that the normalization in this Section differs from the one chosen in Definition~\ref{def_x_i} by
a factor~$\frac{1}{\sqrt{g}}$. The nice feature of this new normalization is that with this choice, the pull-back of divisor classes between some subspaces of the Neron Severi spaces is isometric, see the next Proposition~\ref{prop_phipsi}.
\end{rk}

For~$i\geq 0$ and~$j\geq 1$, one has by Lemma~\ref{Intersections}, with this new normalization for $\gamma^i$,
\begin{align}
\label{norme_gamma_X}&\left\|\gamma_X^i\right\|_X = \sqrt{g_X},
&&\scalaire{\gamma_{X}^i}{\gamma_{X}^{i+j}}_{X} = 
\frac{(q^j + 1) - \sharp \Cun(\FF_{q^j})}{2\sqrt{q^j}},\\
\label{norme_gamma_Y}&\left\|\gamma_Y^i\right\|_Y = \sqrt{g_Y},
&&\scalaire{\gamma_{Y}^i}{\gamma_{Y}^{i+j}}_{Y} = 
\frac{(q^j + 1) - \sharp \Cdeux(\FF_{q^j})}{2\sqrt{q^j}}.
\end{align}
In other words, the vectors $\gamma_{\Cun}^i $ for $i \geq 0$
lie on in the euclidean sphere of radius $\sqrt{g_{\Cun}}$ in the finite dimensional euclidean vector space
$\HVperp_X$, and the data of the scalar products is equivalent to the datas of the numbers of rational points of $\Cun$ on the finite extensions of $\FF_q$.

\begin{prop} \label{prop_phipsi}
Restricted respectively to~$\Frob_X$ and~$\Frob_Y$, the morphisms~$\varphi$ and~$\psi$
satisfy the followings.
\begin{enumerate}
\item One has~$\varphi \circ \psi = \Id_{\Frob_Y}$, the identity map on $\Frob_Y$.
\item \label{phipsi_phigamma}For all $i \geq 0$, $\varphi(\gamma_{\Cun}^i) = \gamma_{\Cdeux}^i$.
\item \label{phipsi_gammadelta} For all
$\gamma \in \Frob_X$ and all $\delta \in \Frob_Y$,
$\scalaire{\gamma}{\psi\left(\delta\right)}_{X}= \scalaire{\varphi(\gamma)}{\delta}_{Y}$ .
\item For all $i, j \geq 0$,
$\scalaire{\psi(\gamma_{\Cdeux}^i)}{\psi(\gamma_{\Cdeux}^j)}_{\Cun}
=
\scalaire{\gamma_{\Cdeux}^i}{\gamma_{\Cdeux}^j}_{\Cdeux}$ .
\item The morphism~$\psi$ is an isometric
embedding~$
\left(\Frob_Y, \scalaire{\cdot}{\cdot}_{\Cdeux}\right)
\hookrightarrow 
\left(\Frob_X, \scalaire{\cdot}{\cdot}_{\Cun}\right)
$
and~$\psi \circ \varphi$ is the orthogonal projection of $\Frob_X$ on its
subspace~$\psi(\Frob_Y)$.
\end{enumerate}
\end{prop}

\begin{preuve} Item (1) follows from the identity 
$(f\times f)_*\circ (f\times f)^* = d^2 \Id_{\Cdeux \times \Cdeux}$. Item (2) follows from the identity $(f\times f)_* \Gamma_{F_{\Cun}^i} = d \Gamma_{F_{\Cdeux}^i}$. Item (3) follows from the projection formula for the proper morphism 
$f\times f : \Cun \times \Cun \rightarrow \Cdeux \times \Cdeux$. Item (4) follows then from items (1) and (3). 
Finally, item (5) follow from the preceding ones together with the following Lemma.
\end{preuve}

\begin{lem}
Let $\left(E,\scalaire{.}{.}_E \right)$ and $\left(F, \scalaire{.}{.}_F \right)$ be two finite dimensional euclidean vector spaces. Suppose that $E$ is generated by a family $\{u_i, i \geq 0\}$ and $F$ is generated by a family $\{v_j, j \geq 0\}$. Let~$\varphi : E \rightarrow F$ and~$\psi : F \rightarrow E$ be linear
maps satisfying~$\varphi \circ \psi = \Id_F$ and such that
\begin{align}
&\label{hyp_isometrie}
\scalaire{\psi(v_i)}{\psi(v_j)}_{E}= \scalaire{v_i}{v_j}_{F},
&& \forall i, j \geq 0,\\
&\label{hyp_adj}
\scalaire{v}{\varphi(u)}_{F}= \scalaire{\psi(v)}{u}_{E},
&&\forall u \in E, \forall v \in F.
\end{align}
Then, $\psi$ is an isometric embedding from $F$ to $E$ and $\psi \circ \varphi$ is the orthogonal projection from $E$ to~$\psi(F)$.
\end{lem}

\begin{preuve}
By~(\ref{hyp_isometrie}), the morphism~$\psi$ must be an isometric embedding.
Let us prove that $\psi \circ \varphi$ is the orthogonal projection
on~$\psi(F)$. Since~$\varphi \circ \psi = \Id_F$, we
have~$(\psi \circ \varphi)^2 = \psi \circ \varphi \circ \psi \circ \varphi = \psi \circ \varphi$; therefore~$\psi \circ \varphi$ is a projector.
For~$v\in F$, one
has~$\psi \circ \varphi(\psi(v)) = \psi(v)$, again
by~$\varphi \circ \psi = \Id_F$; in other terms~$\psi(F)$ is stabilized
by~$\psi\circ\varphi$.
Last, for~$u \in E$, one has~$u-\psi \circ \varphi(u) \in \psi(F)^\perp$.
Indeed, for~$v \in F$, one has:
\begin{align*}
\scalaire{u-\psi \circ \varphi(u)}{\psi(v)}_E 
&=
\scalaire{u}{\psi(v)}_E -\scalaire{\psi \circ \varphi(u)}{\psi(v)}_E\\
&=
\scalaire{\varphi(u)}{v}_F -\scalaire{\varphi(u)}{v}_F
&&\text{by~(\ref{hyp_adj})}\\
&=
0.
\end{align*}
This concludes the proof.
\end{preuve}

\begin{defi}
Let~$\Frob_{X/Y}$ be the subspace~$\psi(\Frob_Y)$, called the {\em relative subspace} for the morphism $X \rightarrow Y$.
\end{defi}

Then we have just proved that
$$
\begin{array}{ccccc}
\Frob_X & = & \Frob_{X/Y} & \oplus & \Frob_{X/Y}^\perp \\
\tikz\node[inner xsep =-2pt,inner ysep = 2pt,rotate=90] {$\in$};
&&
\tikz\node[inner xsep =-2pt,inner ysep = 2pt,rotate=90] {$\in$};
&&
\tikz\node[inner xsep =-2pt,inner ysep = 2pt,rotate=90] {$\in$}; \\
u       & = & \psi \circ \varphi(u) & + & \left(u-\psi \circ \varphi(u)\right)
\end{array}
$$

\begin{rk}
The absolute case of Section~\ref{s_Absolute_bounds} enters in this relative case by choosing any non-constant rational function~$X \rightarrow {\mathbb P}^1$; in this case, the relative space $\Frob_{X/{\mathbb P}^1}$ is the sub-space of $\HVperp_{\Cun}$ generated by the projections of the classes of iterations of the Frobenius morphism. 
\end{rk}

We are now able to easily compute some norms and scalar products.

\begin{lem} \label{normeetpsa2}
For any $i \geq 0$, the following decomposition holds true:
$$
\gamma_{\Cun}^i
=
\underbrace{\psi(\gamma_{\Cdeux}^i)}_{\in \Frob_{X/Y}}
+
\underbrace{\left(\gamma_{\Cun}^i - \psi(\gamma_{\Cdeux}^i)\right)}_{\in \Frob_{X/Y}^\perp}
\in
\Frob_X
$$
and one has:
\begin{align}
&\forall i,j\geq 0,
&&\label{ortho} \psi(\gamma_Y^i) \perp \gamma_X^{j} - \psi(\gamma_Y^{j})\\
&\forall i\geq 0,
&&\label{norme_gamma_X-psi_gamma_Y}\left\|\gamma_{\Cun}^i - \psi(\gamma_{\Cdeux}^i)\right\|_X
=
\sqrt{g_{\Cun}-g_{\Cdeux}} \\
&\forall i\geq 0,\forall j\geq 1,
&&\label{scalaire_gamma_X-psi_gamma_Y}\scalaire{\gamma_{\Cun}^i - \psi(\gamma_{\Cdeux}^i)}{\gamma_{\Cun}^{i+j} - \psi(\gamma_{\Cdeux}^{i+j})}_{\Cun}
= \frac{\sharp \Cdeux(\FF_{q^j})-\sharp \Cun(\FF_{q^j})}{2\sqrt{q^j}}.
\end{align}
\end{lem}

\begin{preuve}
The first norm calculation is just a consequence of Pythagore Theorem:
\begin{align*}
\left\|\gamma_{\Cun}^i\right\|_X^2
&=
\left\|\psi(\gamma_{\Cdeux}^i)\right\|_X^2
+
\left\|\gamma_{\Cun}^i - \psi\left(\gamma_{\Cdeux}^i\right)\right\|_X^2
&& \text{by Pythagore}\\
&=
\left\|\gamma_{\Cdeux}^i\right\|_Y^2
+
\left\|\gamma_{\Cun}^i - \psi\left(\gamma_{\Cdeux}^i\right)\right\|_X^2
&& \text{since $\psi$ isometric}
\end{align*}
and thus by~(\ref{norme_gamma_X}) and~(\ref{norme_gamma_Y})
$$
g_X = g_Y +
\left\|\gamma_{\Cun}^i - \psi\left(\gamma_{\Cdeux}^i\right)\right\|_X^2.
$$
Taking intby o account orthogonality, we also easily compute the scalar product:
\begin{align*}
\scalaire{\gamma_{\Cun}^i - \psi(\gamma_{\Cdeux}^i)}{\gamma_{\Cun}^{i+j} - \psi(\gamma_{\Cdeux}^{i+j})}_{\Cun}
&=
\scalaire{\gamma_X^i}{\gamma_X^{i+j}}_X
-
\scalaire{\psi\left(\gamma_Y^i\right)}{\psi\left(\gamma_Y^{i+j}\right)}_X
&& \text{by orthogonality} \\
&=
\scalaire{\gamma_X^i}{\gamma_X^{i+j}}_X
-
\scalaire{\gamma_Y^i}{\gamma_Y^{i+j}}_Y
&& \text{since $\psi$ isometric}\\
&=
\frac{(q^j + 1) - \sharp \Cun(\FF_{q^j})}{2\sqrt{q^j}}
-
\frac{(q^j + 1) - \sharp \Cdeux(\FF_{q^j})}{2\sqrt{q^j}}
&&\text{by~(\ref{norme_gamma_X}) and~(\ref{norme_gamma_Y})}.
\end{align*}
The result follows.
\end{preuve}

\subsection{Number of points in coverings}

As applications of this Proposition, we prove in the very same spirit than in Section~\ref{s_Absolute_bounds} the following three results. The first one is well known, the others are new. It worth to notice that the first one can be proved using Tate modules of the jacobians of the involved curves (see e.g.~\cite{AubryPerret}), whereas up to our knowledge, the others cannot.

\begin{prop} \label{X-Y}
Suppose that there exists a finite morphism $\Cun \rightarrow \Cdeux$. Then we have
$$\vert \sharp \Cun(\FF_q)-\sharp \Cdeux(\FF_q)\vert \leq 2(g_{\Cun}-g_{\Cdeux})\sqrt q.$$
\end{prop}

\begin{preuve}
We apply Cauchy-Schwartz inequality
to the vectors $\gamma_{\Cun}^0 - \psi(\gamma_{\Cdeux}^0)$
and $\gamma_{\Cun}^1 - \psi(\gamma_{\Cdeux}^1)$.
Thanks to Lemma~\ref{normeetpsa2} specialized to~$j=1$ and $i=0, 1$,
we obtain
\begin{align*}
\left\vert \frac{\sharp \Cun(\FF_{q})-\sharp \Cdeux(\FF_{q})}{2\sqrt{q}}\right\vert^2
&= \left\vert \scalaire{\gamma_{\Cun}^0 - \psi(\gamma_{\Cdeux}^0)}{ \gamma_{\Cun}^1 - \psi(\gamma_{\Cdeux}^1)}_{\Cun} \right\vert^2\\
&\leq
\Vert \gamma_{\Cun}^0 - \psi(\gamma_{\Cdeux}^0) \Vert_{\Cun}^2 \times \Vert \gamma_{\Cun}^1 - \psi(\gamma_{\Cdeux}^1)\Vert_{\Cun}^2\\
&=
(g_{\Cun}-g_{\Cdeux})^2,
\end{align*}
hence the Proposition.
\end{preuve}

%
%

The following Proposition is the relative form of Proposition~\ref{X(Fq2)}. 
Of course, although less nice, such upper bounds can be given for any $n$.

\begin{prop} \label{RelatifNonTrivial}
For any finite morphism $\Cun \rightarrow \Cdeux$ with $g_{\Cun}\neq g_{\Cdeux}$, we have
$$\sharp \Cun({\mathbb F}_{q^2})- \sharp \Cdeux({\mathbb F}_{q^2}) \leq 2(g_{\Cun}-g_{\Cdeux}) q - \frac{\Bigl(\sharp \Cun({\mathbb F}_{q})-\sharp \Cdeux({\mathbb F}_{q})\Bigr)^2}{g_{\Cun}-g_{\Cdeux}}.$$
\end{prop}

\begin{preuve}
After computating the $3\times 3$ Gram determinant 
$$G_3(X/Y) = \Gram(\gamma_X^{0} - \psi(\gamma_Y^{0}), \gamma_X^{1} - \psi(\gamma_Y^{1}), \gamma_X^{2} - \psi(\gamma_Y^{2}))$$
 using Lemma~\ref{normeetpsa2}, it is readily seen, with the very same proof, that factorization in Lemma~\ref{lem_factorization} holds also in this relative case, so that there exists two (explicitly given as determinants) factors $G_{n}^+(\Cun/\Cdeux)$ and $G_{n}^-(\Cun/\Cdeux)$ ($n=2$ is sufficient here) as in Definition~\ref{def_minors}, such that $G_{n}(\Cun/\Cdeux)= G_{n}^+(\Cun/\Cdeux) \times G_{n}^-(\Cun/\Cdeux)$  also holds. 
 Moreover, Proposition~\ref{DnRecursif}
 also continue to hold, hence for $n=2$, the result follows from
$G_2^-(\Cun/\Cdeux) \geq 0$.
\end{preuve}

\subsection{Number of points in a fiber product}
Let
\begin{equation} \label{DiagrXY1Y2Z}
\begin{tikzpicture}[>=latex,baseline=(M.center)]
\matrix (M) [matrix of math nodes,row sep=0.65cm,column sep=0.75cm]
{
                              &  |(Y1)| Y_1 & \\
|(Y1xY2)| X  &            & |(Z)| Z \\
                              &  |(Y2)| Y_2 & \\
};
\draw[->] (Y1xY2) -- (Y1) node[midway,anchor = south east] {$p_1$} ;
\draw[->] (Y1) -- (Z) node[midway,anchor = south west] {$f_1$} ;
\draw[->] (Y1xY2) -- (Y2) node[midway,below,anchor = north east] {$p_2$} ;
\draw[->] (Y2) -- (Z) node[midway,anchor = north west] {$f_2$} ;
\end{tikzpicture}
\end{equation}
be a commutative diagram of finite covers
of absolutely irreducible smooth projective curves defined over~$\FF_q$.
Applying the results of the beginning of this Section to the four morphisms
involved
in this diagram leads to another diagram between Euclidean spaces
\begin{equation*}\label{XY1Y2Z}
\begin{tikzpicture}[baseline=(FZ),node distance = 2.5cm, bend angle=30]
\node (FY1) {$\Frob_{Y_1}$};
\node[below left of = FY1] (FZ) {$\Frob_Z$};
\node[below right of = FZ] (FY2) {$\Frob_{Y_2}$};
\node[below right of = FY1] (FX) {$\Frob_X$};
\draw[thick,every node/.style={font=\sffamily\tiny}]
  (FZ) edge[right hook-latex] node[midway,below,sloped] {$\psi_{Y_1/Z}$} (FY1)
  (FZ) edge[right hook-latex] node[midway,above,sloped] {$\psi_{Y_2/Z}$} (FY2)
  (FY1) edge[right hook-latex] node[midway,below,sloped] {$\psi_{X/Y_1}$} (FX)
  (FY2) edge[right hook-latex] node[midway,above,sloped] {$\psi_{X/Y_2}$} (FX)
  (FY1.west) edge[dashed,bend right,->,>=latex] node[midway,above,sloped] {$\varphi_{Y_1/Z}$} (FZ.north) 
  (FY2.west) edge[dashed,bend left,->,>=latex] node[midway,below,sloped] {$\varphi_{Y_2/Z}$} (FZ.south)
  (FX.north) edge[dashed,bend right,->,>=latex] node[midway,above,sloped] {$\varphi_{X/Y_1}$} (FY1.east)
  (FX.south) edge[dashed,bend left,->,>=latex] node[midway,below,sloped] {$\varphi_{X/Y_2}$} (FY2.east) ;
\end{tikzpicture}
\end{equation*}
\medskip
Let us introduce the following hypothesis:
$$\hbox{\it The fiber product $Y_1\times_Z Y_2$ is abolutely irreducible and smooth.}\leqno(H)$$
\medskip

\begin{prop} \label{prop_XY12Z} Suppose that $(H)$ holds, and that $X = Y_1\times_Z Y_2$ in diagram~\eqref{DiagrXY1Y2Z}.
Then
\begin{equation*}\label{phi_psi_XY12Z}
\varphi_{X/Y_2} \circ \psi_{X/Y_1}
=
\psi_{Y_2/Z} \circ \varphi_{Y_1/Z} \\
\end{equation*}
on~$\Frob_{Y_1}$.
\end{prop}

\begin{preuve}
The proof of formula~(\ref{phi_psi_XY12Z}) is mainly set theoretic.
We write: 
\begin{align*}
(p_2\times p_2)_*\circ (p_1\times p_1)^*\left(\Gamma_{F_{\Cdeux_1}}\right)
&=
(p_2\times p_2)_*\left(\{ [(y_1, y_2), (y'_1, y'_2)] \in (\Cdeux_1\times_Z\Cdeux_2)^2 ; y'_1=F_{\Cdeux_1}(y_1)\}\right)\\
&= 
\{ (y_2, y'_2) \in \Cdeux_2\times \Cdeux_2 ; \exists y_1 \in \Cdeux_1 ; f_1(y_1)=f_2(y_2) ; f_2(y'_2)=F_{\Cdeux_2}(f_1(y_1))\}\\
&= \deg f_1 \times \{ (y_2, y'_2) \in \Cdeux_2\times \Cdeux_2 ;  f_2(y'_2)=F_{\Cdeux_2}(f_2(y_2))\}\\
&= \deg f_1 \times (f_2\times f_2)^*(\Gamma_{F_{Z}}).
\end{align*}
Now, since $Y_1\times_{\Cun} Y_2$ is assumed to be absolutely irreducible, we have $\deg f_1=\deg p_2$ and $\deg f_2 = \deg p_1$ in a fiber product setting. Taking into account  the normalization (\ref{phipsi}), and projecting
onto~$\HVperp_{Y_2}$ as in Section 4,  allow us to conclude.
\end{preuve}

This  allow us us to compute some other norms and scalar products.

\begin{lem} \label{lem_orthogonalite_XY12Z}
Suppose that $(H)$ holds, and that $X = Y_1\times_Z Y_2$ in diagram~\eqref{DiagrXY1Y2Z}. Then we have in~$\Frob_X$
\begin{multline*}
\psi_{X/Y_1}(\gamma_{Y_1}^i),\psi_{X/Y_2}(\gamma_{Y_2}^i),\psi_{X/Z}(\gamma_{Z}^i)\\
\in
\Vect\left(
\psi_{X/Y_{1}}\left(\gamma_{Y_1}^i-\psi_{Y_1/Z}(\gamma_Z^i)\right),
\psi_{X/Y_{2}}\left(\gamma_{Y_2^i}-\psi_{Y_2/Z}(\gamma_Z^i)\right),
\gamma_{X}^i-\psi_{X/Z}(\gamma_Z^i)
\right)^\perp
\end{multline*}
and
\begin{align*}\label{phi_psi_XY12Z_2}
&\forall i, j \in \NN,
&&
\psi_{\Cun/\Cdeux_1} \left(\gamma_{\Cdeux_1}^i- \psi_{\Cdeux_1/Z}(\gamma_Z^i)\right)
\perp
\psi_{\Cun/\Cdeux_2} \left(\gamma_{\Cdeux_2}^j- \psi_{\Cdeux_2/Z}(\gamma_Z^j)\right).
\end{align*}
\end{lem}

\begin{preuve}
First, we compute scalar products with~$\psi_{X/Y_1}(\gamma_{Y_1}^i)$. We have
\begin{align*}
\scalaire
{\psi_{X/Y_1}(\gamma_{Y_1}^i)}
{\psi_{X/Y_{1}}\left(\gamma_{Y_1}^i-\psi_{Y_1/Z}(\gamma_Z^i)\right)}_X
&=
\scalaire{\gamma_{Y_1}^i}
{\gamma_{Y_1}^i-\psi_{Y_1/Z}(\gamma_Z^i)}_{Y_1}
&& \text{since $\psi_{X/Y_1}$ is isometric} \\
&= 0,
&& \text{by Lemma~\ref{normeetpsa2}}
\end{align*}
\begin{align*}
\scalaire{\psi_{X/Y_1}(\gamma_{Y_1}^i)}
{\psi_{X/Y_{2}}\left(\gamma_{Y_2^i}-\psi_{Y_2/Z}(\gamma_Z^i)\right)}_X
&=
\scalaire{\varphi_{X/Y_2}\circ\psi_{X/Y_1}(\gamma_{Y_1}^i)}
{\gamma_{Y_2^i}-\psi_{Y_2/Z}(\gamma_Z^i)}_{Y_2}
&&\text{by prop.~\ref{prop_phipsi}}\\
&=
\scalaire{\psi_{Y_2/Z}\circ\varphi_{Y_1/Z}(\gamma_{Y_1}^i)}
{\gamma_{Y_2^i}-\psi_{Y_2/Z}(\gamma_Z^i)}_{Y_2}
&&\text{by prop.~\ref{prop_XY12Z}}\\
&=
\scalaire{\psi_{Y_2/Z}(\gamma_{Z}^i)}
{\gamma_{Y_2^i}-\psi_{Y_2/Z}(\gamma_Z^i)}_{Y_2}
&&\text{by prop.~\ref{prop_phipsi}}\\
&=
0,
&& \text{by Lemma~\ref{normeetpsa2}}
\end{align*}
and
\begin{align*}
\scalaire
{\psi_{X/Y_1}(\gamma_{Y_1}^i)}
{\gamma_{X}^i-\psi_{X/Z}(\gamma_Z^i)}_X
&=
\scalaire
{
\psi_{X/Y_1}(\gamma_{Y_1}^i)
}
{
\gamma_{X}^i-\psi_{X/Y_1}(\gamma_{Y_1}^i)
+
\psi_{X/Y_1}\left(\gamma_{Y_1}^i - \psi_{Y_1/Z}(\gamma_Z^i)\right)
}_X \\
&= 0.
\end{align*}
The same holds for~$\psi_{X/Y_2}(\gamma_{Y_2}^i)$.

Secondly, we compute
scalar products
with~$\psi_{X/Z}(\gamma_{Z}^i) = \psi_{X/Y_\epsilon}\circ\psi_{Y_\epsilon/Z}(\gamma_Z^i)$,
for~$\epsilon=1,2$. We have:
\begin{align*}
\scalaire{\psi_{X/Z}(\gamma_{Z}^i)}
{\psi_{X/Y_{\epsilon}}\left(\gamma_{Y_\epsilon}^i-\psi_{Y_\epsilon/Z}(\gamma_Z^i)\right)}_X
&=
\scalaire{\psi_{X/Y_\epsilon}\circ\psi_{Y_\epsilon/Z}(\gamma_Z^i)}
{\psi_{X/Y_{\epsilon}}\left(\gamma_{Y_\epsilon}^i-\psi_{Y_\epsilon/Z}(\gamma_Z^i)\right)}_X \\
&=
\scalaire{\psi_{Y_\epsilon/Z}(\gamma_Z^i)}
{\gamma_{Y_\epsilon}^i-\psi_{Y_\epsilon/Z}(\gamma_Z^i)}_{Y_\epsilon}
&&\text{$\psi_{X/Y_\epsilon}$ isometric} \\
&= 0.
&&\text{by Lemma~\ref{normeetpsa2}}
\end{align*}
Last,~$\scalaire
{\psi_{X/Z}(\gamma_{Z}^i)}
{\gamma_{X}^i-\psi_{X/Z}(\gamma_Z^i)}_X
= 0$ by Lemma~\ref{normeetpsa2}.
\end{preuve}

In the same way than the use of the orthogonal projections associated
to~$\Frob_{X/Y} \subset \Frob_X$ in the preceding Section,
we do the same here with~$\Frob_{X/Y_1} + \Frob_{X/Y_2} \subset \Frob_X$.

\begin{lem}
Suppose that $(H)$ holds, and that $X = Y_1\times_Z Y_2$ in diagram~\eqref{DiagrXY1Y2Z}. For~$i\geq 0$, let~$\gammaa^i$ be the orthogonal projection
of~$\gamma_{X}^i$ onto~$\left(\Frob_{X/Y_1}+\Frob_{X/Y_2}\right)^\perp$
inside~$\Frob_X$. Then
\begin{equation*}\label{proja4}
\gammaa^i
=
\gamma^i_{X} - \psi_{X/\Cdeux_1}(\gamma^i_{\Cdeux_1}) - \psi_{X/\Cdeux_2}(\gamma^i_{\Cdeux_2}) + \psi_{X/Z}(\gamma^i_{Z}),
\end{equation*}
and one has:
\begin{align*}
\Vert \gammaa^i\Vert_X
&=
\sqrt{g_{X}-g_{\Cdeux_1}-g_{\Cdeux_2}+g_Z}\\
\scalaire{\gammaa^i}{\gammaa^{i+j} }_{X}
&= \frac{\sharp Y_1(\FF_{q^j})+\sharp Y_2(\FF_{q^j})- \sharp X(\FF_{q^j})- \sharp Z(\FF_{q^j})}{2\sqrt{q^j}}.
\end{align*}
\end{lem}

\begin{preuve}
Clearly~$\gamma_X^i - \gammaa^i \in \Frob_{X/Y_1} + \Frob_{X/Y_2}$.
Thanks to orthogonality relations of Lemma~\ref{lem_orthogonalite_XY12Z}
and rewriting~$\gammaa^i$ a little bit like
$$
\gammaa^i
=
\left(\gamma^i_{X} - \psi_{X/Z}(\gamma_Z^i)\right)
- \psi_{X/\Cdeux_1}\left(\gamma^i_{\Cdeux_1}-\psi_{Y_1/Z}(\gamma^i_{Z})\right)
- \psi_{X/\Cdeux_2}\left(\gamma^i_{\Cdeux_2}) -\psi_{Y_2/Z}(\gamma^i_{Z})\right),
$$
we prove that~$\gammaa^i \in \left(\Frob_{X/Y_1} + \Frob_{X/Y_2}\right)^\perp$.

To compute the norm let us note that in the decomposition
$$
\gamma^i_{X} - \psi_{X/Z}(\gamma_Z^i)
=
\left[
\psi_{X/\Cdeux_1}\left(\gamma^i_{\Cdeux_1}-\psi_{Y_1/Z}(\gamma^i_{Z})\right)+
\psi_{X/\Cdeux_2}\left(\gamma^i_{\Cdeux_2}) -\psi_{Y_2/Z}(\gamma^i_{Z})\right)
\right]
+
\left[\gammaa^i\right].
$$
The first bracket lies on in~$\Frob_{X/Y_1} + \Frob_{X/Y_2}$, while the second
one lies on in~$\left(\Frob_{X/Y_1} + \Frob_{X/Y_2}\right)^\perp$. Applying Pythagore
again, we deduce that
\begin{align*}
\left\|\gamma^i_{X} - \psi_{X/Z}(\gamma_Z^i)\right\|_X^2
&=
\left\|
\psi_{X/\Cdeux_1}\left(\gamma^i_{\Cdeux_1}-\psi_{Y_1/Z}(\gamma^i_{Z})\right)+
\psi_{X/\Cdeux_2}\left(\gamma^i_{\Cdeux_2} -\psi_{Y_2/Z}(\gamma^i_{Z})\right)
\right\|_X^2
+
\left\|
\gammaa^i
\right\|_X^2 \\
&=
\left\|
\psi_{X/\Cdeux_1}\left(\gamma^i_{\Cdeux_1}-\psi_{Y_1/Z}(\gamma^i_{Z})\right)
\right\|_X^2
+
\left\|
\psi_{X/\Cdeux_2}\left(\gamma^i_{\Cdeux_2}) -\psi_{Y_2/Z}(\gamma^i_{Z})\right)
\right\|_X^2
+
\left\|
\gammaa^i
\right\|_X^2 \\
&=
\left\|
\gamma^i_{\Cdeux_1}-\psi_{Y_1/Z}(\gamma^i_{Z})
\right\|_{Y_1}^2
+
\left\|
\gamma^i_{\Cdeux_2} -\psi_{Y_2/Z}(\gamma^i_{Z})
\right\|_{Y_2}^2
+
\left\|
\gammaa^i
\right\|_X^2.
\end{align*}
Lemma~\ref{normeetpsa2}  allows us to conclude since
\begin{align*}
\left\|
\gammaa^i
\right\|_X^2
&=
\left\|\gamma^i_{X} - \psi_{X/Z}(\gamma_Z^i)\right\|_X^2
-
\left\|
\gamma^i_{\Cdeux_1}-\psi_{Y_1/Z}(\gamma^i_{Z})
\right\|_{Y_1}^2
-
\left\|
\gamma^i_{\Cdeux_2} -\psi_{Y_2/Z}(\gamma^i_{Z})
\right\|_{Y_2}^2 \\
&=
\left(g_X - g_Z\right) - \left(g_{Y_1} - g_Z\right) - \left(g_{Y_2} - g_Z\right).
\end{align*}

In the same way, the calculation of~$\scalaire{\gammaa^i}{\gammaa^{i+j}}$
is a consequence
of the computation
of
the scalar product~$\scalaire{\gamma^i_{X} - \psi_{X/Z}(\gamma_Z^i)}{\gamma^{i+j}_{X} - \psi_{X/Z}(\gamma_Z^{i+j})}$.
\end{preuve}

Last we can prove the following results.

\begin{theo} \label{final}
Let $X, Y_1, Y_2$ and $Z$ be absolutely irreducible smooth projective curves in a cartesian diagram~\eqref{DiagrXY1Y2Z} of finite morphisms. Suppose that the fiber
product~$Y_1\times_Z Y_2$ is also absolutely irreducible and smooth. Then
$$
\left\vert \sharp \Cun(\FF_q)-\sharp \Cdeux_1(\FF_q)-\sharp \Cdeux_2(\FF_q)+\sharp Z(\FF_q)\right\vert
\leq 2(g_{\Cun}-g_{\Cdeux_1}-g_{\Cdeux_2}+g_Z) \sqrt q.
$$
\end{theo}

\begin{preuve} First, for~$X = Y_1\times_Z Y_2$ the result is a direct
consequence of Cauchy-Schwartz inequality applied to $\gammaa^0$
and~$\gammaa^1$ thanks to Lemma~\ref{lem_orthogonalite_XY12Z} since $(H)$ holds by assumption.

Now, for general $\Cun$ satisfying the assumptions of the Theorem, we have by the universal property of the fibered product a finite morphism $\Cun \rightarrow Y_1\times_Z Y_2$. By triangular inequality, and using Proposition~\ref{X-Y} together with the result for $X$, we have
\begin{align*}
\left\vert \sharp \Cun(\FF_q)-\sharp \Cdeux_1(\FF_q)-\sharp \Cdeux_2(\FF_q)+\sharp Z(\FF_q)\right\vert
& \leq
\vert \sharp \Cun({\mathbb F}_q) - \sharp Y_1\times_Z Y_2({\mathbb F}_q)\vert \\
&+ \left\vert \sharp Y_1\times_Z Y_2(\FF_q)-\sharp \Cdeux_1(\FF_q)-\sharp \Cdeux_2(\FF_q)+\sharp Z(\FF_q)\right\vert\\
&  \leq 2(g_{\Cun}-g_{Y_1\times_Z Y_2}) \sqrt q + 2(g_{Y_1\times_Z Y_2}-g_{\Cdeux_1}-g_{\Cdeux_2}+g_Z) \sqrt q\\
&= 2(g_{\Cun}-g_{\Cdeux_1}-g_{\Cdeux_2}+g_Z) \sqrt q
\end{align*} 
and the Theorem is proved.
\end{preuve}

\medskip

It worth to notice that Theorem~\ref{final} cannot holds without any hypothesis. For instance, if $X=Y_1=Y_2$ and the morphisms $Y_i \rightarrow Z$ are equal, then the right hand side equals $2(g_X-2g_X+g_Z)\sqrt{q}=-2(g_X-g_Z)\sqrt{q}$, a negative number! In this case, the Theorem doesn't apply since  $Y_1\times_Z Y_2$ is not irreducible.

\begin{rk}
For~$Y_1 \times_Z Y_2$ to be absolutely irreducible, it suffices that
the tensor product of function fields~$\FF_q(Y_1) \otimes_{\FF_q(Z)} \FF_q(Y_2)$
to be an integral domain. For~$(Q_1,Q_2) \in Y_1 \times_Z Y_2$ to be smooth
it is necessary and sufficient that at least one of the
morphism~$Y_i \to Z$ is unramified at~$Q_i$.
\end{rk}

\section{Questions}\label{s_questions}
In this Section, we suggest a few questions raised by the viewpoint of this article.

\begin{enumerate}
\item Numerical calculations using our algorithm make it possible that the genus sequence $(g_n, n\in \NN)$ whose existence is asserted in Proposition~\ref{Suite_g_n} has, at least asymptotically, a quite nice behaviour. It seems that for fixed $n \geq 1$, there do exists a constant $0<c_n<1$, such that for $q$ large,
$$g_n \sim c_n \sqrt{q}^n,$$
with moreover $\lim_{n \to +\infty} c_n=1$. Of course, we have proved that $c_2= \frac{1}{2}$ and $c_3= \frac{1}{\sqrt{2}}$.

It numerically also seems  that for any $n \geq 1$ and any $q$,
$$g_n < \frac{\sqrt{q}^{n+1}}{\sqrt{q}-1}.$$
\item Is it true that the first requirement of criteria~\ref{criteria_min_x_1} implies the second one, in the same way that it do implies the third one as proved in corollary~\ref{DeriveeDirectionnelle} ? It turns out that numerically, our algorithm shows that this holds true for all pairs $(g, q)$ we have tested. We gave up this question when we became aware that our bounds seems to be equivalent to Oesterl\'e's.
\item Why do the upper bounds obtained in this article always coincide numerically with Osterl\'e's one ?
%
%
\item In view of the hope to write down explicitly Weil bounds of order $n \geq 4$, it is necessary to solve a one variable polynomial equation on $x_1$ of degree $\lceil \frac{n+1}{2}\rceil$, which is greater than 5 for $n \geq 8$. This polynomial equation is, eliminating $\frac{1}{g}$, the resultant
$$\hbox{Res}_{\frac{1}{g}}\left(G_n^-(P_n^g(x_1)), G_n^+(P_n^g(x_1))\right) \in \QQ\left(\frac{1}{{\sqrt q}}\right)[x_1].$$
Do this polynomial have solvable Galois group over the rational function field $\QQ\left(\frac{1}{{\sqrt q}}\right)$ ?
\item Is there a nice interpretation in this framework for the number of rational points of the Jacobian $\hbox{Jac}(X)$ of a curve $X$ ? In the affirmative, same question for the Prym variety associated to an unramified double cover $X \rightarrow Y$ ?
\item Do the viewpoint of this article can be extended to higher dimensional varieties ?
\item As explained in the article, the geometry and the arithmetic of $X$ are traduced via Hodge Index Theorem to an euclidean property throught Gram determinants by $P_n(X) \in \overline{\Dcal}_n$, and its arithmetic by 
$P_n(X) \in \Hcal_n^g$ through Ihara's constraints, so that upper bounds for $\sharp X(\FF_q)$ are derived from the lower bound
for $x_1$ on $\overline{\Dcal}_n \cap \Hcal_n^g$. 

Any other constraint coming either from the geometry or the arithmetic of $X$ is likely to reduce the domain $\Dcal_n \cap \Hcal_n^g$, making possible that the $x_1$ function on the new domain as greater lower bound.
For instance, 
do the records obtained in the table {\tt http://www.manypoints.org/}, better than Osterl\'e bound, can be understood in this way ?
\end{enumerate}

\providecommand{\bysame}{\leavevmode\hbox to3em{\hrulefill}\thinspace}
\providecommand{\MR}{\relax\ifhmode\unskip\space\fi MR }
\providecommand{\MRhref}[2]{%
  \href{http://www.ams.org/mathscinet-getitem?mr=#1}{#2}
}
\providecommand{\href}[2]{#2}

\end{document}